\documentclass[11pt,amssymb]{article}
\usepackage{amsfonts,latexsym}

\newcommand{\R}{\mathbb R}
\newcommand{\Z}{\mathbb Z}
\newcommand{\N}{\mathbb N}

\newcommand{\LL}{\mathbb L}

\newtheorem{theorem}{Theorem}
\newtheorem{lemma}{Lemma}
\newtheorem{proposition}{Proposition}
\newtheorem{corollary}{Corollary}
\newtheorem{definition}{Definition}

\def\ep{\varepsilon}

\begin{document}

\title{ Uniqueness of positive periodic solutions with some peaks.}
\author{Genevi\`eve Allain $^a$ and Anne Beaulieu $^b$.}

\maketitle

   $^a$ Laboratoire d'Analyse et de Math\'ematiques Appliqu\'ees, Facult\'e de Sciences et Technologie,
                                Universit\'e Paris-Est Cr\'eteil, Cr\'eteil, France.\\
            genevieve.allain@u-pec.fr\\                   
       
  $^b$ Laboratoire d'Analyse et de Math\'ematiques Appliqu\'ees, Universit\'e Paris-Est Marne la Vall\'ee,
   Marne la Vall\'ee, France.\\
 anne.beaulieu@univ-mlv.fr\\

{\bf Abstract.} This work deals with the semilinear equation $-\Delta u+u-u^p=0$ in $\R^N$, $2\leq p<{N+2\over N-2}$. We consider the positive solutions which are ${2\pi\over\ep}$-periodic in $x_1$ and decreasing to 0 in the other variables, uniformly in $x_1$. 
Let a periodic configuration of points be given on the $x_1$-axis, which repel each other as the period tends to infinity. If there exists a solution which has these points as peaks, we prove that the points must be asymptotically uniformly distributed on the $x_1$-axis. Then,
 for $\ep$ small enough, we prove the uniqueness up to a translation of the positive solution with some peaks on the $x_1$-axis, for a given minimal period in $x_1$, and we estimate the difference between this solution and the groundstate solution.\\

 

\section{Introduction.}
We consider the equation
\begin{equation}\label{eq:u}
-\Delta u+u-u_+^p=0\hbox{   in   } \frac{S^1}{\ep}\times\R^{N-1}\\
\end{equation}
where $u_+=\max(u,0).$\\
By $\frac{S^1}{\ep}$, we mean that 
$$u(x_1+\frac{2\pi}{\ep},x')=u(x_1,x')$$ and that 
$${\partial u\over\partial x_1}(x_1+\frac{2\pi}{\ep},x')={\partial u\over\partial x_1}(x_1,x').$$
We suppose that 
$$u(x_1,x')\rightarrow0,\quad\hbox{ as $\vert x'\vert\rightarrow+\infty$, uniformly in $x_1$}. $$
If $u>0$, then $u$ is radial and decreasing in $x'$. Let us indicate how the moving plane method (\cite{Gidas}, \cite{Berestycki}) can be applied to prove this property.\\

First, following the proof of Proposition 1.4 in \cite{Gidas}, we define
$$\overline h(x')=\int_{S^1\over\ep}\int_{S^{N-2}}u(x_1,\vert x'\vert\theta)d\theta dx_1$$
and we obtain
$$\overline h(x')\leq Ce^{-a\vert x'\vert},$$
for some $a>0$ and some $C$ independent of $x$. Then the Harnack inequalities (see \cite{GT}) give
\begin{equation}\label{eq:cptexp}
u(x_1,x')\leq Ce^{-a\vert x'\vert},
\end{equation}
for some $C$ independent of $x$.\\
Now, for all $i=2,...,N$, we denote $x=(x_i,y)$. For a given $t\in\R$, we define
$$Q_t=\{x\in \frac{S^1}{\ep}\times\R^{N-1};x_{i}<
t\};
\quad u_t(x)=u(2t-x_i,y)$$
and
$$\Lambda=\{t\in\R;\forall\mu>t,u\geq u_{\mu}\quad\hbox{in $Q_{\mu}$}\}.$$
The first step of the moving plane method is to prove that $\Lambda$ is nonempty and bounded from below. The second step is to show that if $t_0=\hbox{inf $\Lambda$}$, then $u=u_{t_{0}}$ in $Q_{t_0}$. 
The proof of the first step is standard. We multiply 
$$-\Delta(u_t-u)=-u_t+u+u_t^p-u^p$$
by $(u_t-u)+$ and we integrate on $Q_t$.
We prove the second step as follows.\\
By continuity, $u\geq u_{t_0}$ in $Q_{t_0}$. Then the strong maximum principle yields either $u=u_{t_0}$  in $Q_{t_0}$ or $u>u_{t_0}$ in $Q_{t_0}$.
If
 $$f(u)=-u+u^p.$$
For all $\eta>0$, we have
 $$\vert f(u_t)-f(u)\vert(u_t-u-\eta)_+^2\leq C(u_t+u_t^p)(u_t-u-\eta)_+^2.$$
 with $C$ independent of $\eta$.\\
 
 Moreover, we know by (\ref{eq:cptexp}) that $u\in L^r$ for all $1\leq r\leq\infty$. Consequently, we can argue as in \cite{Ge}, section 2.2. With the notations there, we let $q=2$ and $z=(u_t-u-\eta)_+^2$. We get
 $$
  \int_{Q_t\backslash K}\vert \nabla z\vert^2dx\leq C(\int_{Q_t\backslash K}(u_t+u_t^p)z^2dx
  $$
where $K$ is any compact set in $Q_t$ such that $(u_t-u)_+=0$ in $K$ and $C$ is independent of $t$, $\eta$ and of $K$.
If $\delta>0$ and $T>0$ are given, we define
 $$ K'=\{x'\in\R^{N-1};-T\leq x_j\leq t_0-\delta;\quad j=2,...,N\}$$
 and 
 $$K=\frac{S^1}{\ep}\times K'.$$
  Using H\"{o}lder inequality, we obtain a real number $C$ independent of $t$, $\delta$ and $T$ such that
 $$
 \int_{Q_t\backslash K}\vert \nabla z\vert^2dx\leq C(\int_{S^1\over\ep}(\int_{ K'}(u_t+u_t^p)^{N-1\over2}dx')^{2\over N-1}(\int_{K'}z^{2^{\star}}dx')^{2\over 2^{\star}}dx_1
 $$
 where $2^{\star}={2(N-1)\over N-3}$.\\
 Now, since $z$ has compact support, we can use the Sobolev-Gagliardo-Nirenberg inequality to get for all $x_1$
 $$(\int_{K'}z^{2^{\star}}(x_1,x')dx')^{2\over 2^{\star}}\leq C\int_{K'}\vert\nabla'z\vert^2(x_1,x') dx',$$
 where $C$ is independent of $x_1$, $t$, $\delta$, $\eta$ and $T$.
 We are led to

 $$
 \int_{Q_t\backslash K}\vert \nabla z\vert^2dx\leq C\int_{S^1\over\ep}(\int_{ K'}(u_t+u_t^p)^{N-1\over2}dx')^{2\over N-1}\int_{K'}\vert\nabla z\vert^2(x_1,x') dx')dx_1
$$
 
 On one hand, we choose $\delta$ small enough and $T$ large enough to have for all $t_0-\delta<t<t_0$.
 $$
\forall x_1,\quad C\int_{ K'}(u_t+u_t^p)^{N-1\over2}dx')^{2\over N-1}dx_1<\frac{1}{2}.
$$
Letting $\eta\rightarrow0$, we get
\begin{equation}\label{eq:inge2}
\int_{Q_t\backslash K}\vert \nabla(u_t-u)_+^2\vert^2dx\leq {1\over2}\int_{Q_t\backslash K}\vert \nabla(u_t-u)_+^2\vert^2.
\end{equation}

On the other hand, since $K$ is a compact set and since $K\subset Q_{t_0}$,  then, if $u>u_{t_0}$ in $Q_{t_0}$, there exists $0<\delta_1<\delta$ such that
 $u>u_t\hbox{ in $K$ for all $t\in]t_0-\delta_1,t_0[$}.$
  In view of (\ref{eq:inge}) and (\ref{eq:inge2}), we deduce that the case $u>u_{t_0}$  in $Q_{t_0}$ leads to $u\geq u_t$ in $Q_t$ for all $t\in]t_0-\delta_1,t_0[$ and this is in contradiction with the definition of $t_0$. Consequently, $u=u_{t_0}$ in $Q_{t_0}$. The proof that $u$ is radial and decreasing in $x'$ is complete.\\
  
 We consider the subcritical case
$$2\leq p<\frac{N+2}{N-2}\quad\hbox{for $N\geq3$,}\quad p\geq 2\quad\hbox{for $N=2$}.$$
We assume that $p\geq 2$ instead of $p>1$ for some technical reasons, appearing in section 3.\\
Let $U$ be the groundstate solution in $\R^N$. It verifies
$$-\Delta U+U-U^p=0\quad\hbox{in $\R^N$}.$$
It is known that $U$ is positive, radial and decreasing to 0 at infinity. Moreover the behavior at infinity is
$$\vert x\vert^{\frac{N-1}{2}}e^{\vert x\vert}U(x)\rightarrow L_0\quad\hbox{as $\vert x\vert\rightarrow+\infty$}$$
and
$$\vert x\vert^{\frac{N-1}{2}}e^{\vert x\vert}\frac{\partial U}{\partial x_1}(x)\rightarrow L_1\quad\hbox{as $\vert x\vert\rightarrow+\infty$, $x_1>0$,}$$
for some positive limits $L_0$ and $L_1$.
 (see \cite{Kwong}.)\\
Several recent articles deal with the construction of positive solutions for the equation
$$-\Delta u+u-u^p=0\quad\hbox{in $\R^N$.}$$
Let us refer to 
\cite{Dancer}, \cite{Malchiodi}, \cite{delpacard}.\\
Let us call the Dancer solution the positive solution of (\ref{eq:u}) which is $\frac{2\pi}{\ep}$-periodic in $x_1$, tending to 0 as $\vert x'\vert\rightarrow+\infty$, even in $x_1$ and decreasing in $x_1$ in $[0,\frac{\pi}{\ep}]$. This solution, that we call $u_D$ was constructed in \cite{Dancer} by a bifurcation from the ground-state solution in $\R^{N-1}$. The Dancer solution exists when $0<\ep<\ep^{\star}$, where $\ep^{\star}$ is a known threshold. We have  
 \begin{equation}\label{eq:Dancer}
 \|u_D-U\|_{L^{\infty}(]-\frac{\pi}{\ep},\frac{\pi}{\ep}[\times\R^{N-1})}\rightarrow0\hbox{  as $\ep\rightarrow0$}.
 \end{equation}
 For all $x'$ the fonction $x_1\mapsto u_D(x_1,x')$ reaches its maximum value at the points $\frac{2l\pi}{\ep}$ and reaches its minimum value at the points $\frac{(2l+1)\pi}{\ep}$, $l\in\Z$.

Now, for any $\ep>0$, for any $k\geq2$, let $a_{\ep}^i$, $i=1$,...,$k$, be $k$ points of $[-\pi,\pi[$
which are such that
\begin{equation}\label{eq:hyp}
\frac{ a_{\ep}^{i+1}- a_{\ep}^{i}}{\ep}\rightarrow+\infty
\quad\hbox{as $\ep\rightarrow0$, $i=0,...,k$}.
\end{equation}
where we denote $a_{\ep}^0=a_{\ep}^k-2\pi$ and $a_{\ep}^{k+1}=a_{\ep}^0+2\pi$. \\
Let us denote
$$U_i(x_1,x')=U(x_1-{a_{\ep}^i\over\ep},x').$$
Let us give the following

\begin{definition}\label{pics}
The solution $u$ of (\ref{eq:u}) admits the points $\frac{a_{\ep}^1}{\ep}$,...,$\frac{a_{\ep}^k}{\ep}$ as peaks if 
$\frac{a_{\ep}^1}{\ep}$,...,$\frac{a_{\ep}^k}{\ep}$ are $k$ points of $[-\frac{\pi}{\ep},\frac{\pi}{\ep}[$ 
 verifying (\ref{eq:hyp}) and if
  \begin{equation}\label{eq:conv}
\| u- \sum_{i=1}^k U_i\|_{L^{\infty}(]-\frac{\pi}{\ep},\frac{\pi}{\ep}[\times\R^{N-1})}\rightarrow0 \quad\hbox{as $\ep\rightarrow0$}.
\end{equation}
\end{definition}
Let us remark that by the Maximum Principle, any solution of (\ref{eq:u}) verifying (\ref{eq:conv}) needs to be positive.\\
We can ask whether for any configuration of points in a period which repel each other in the sense of (\ref{eq:hyp}), there exists a solution having these points as peaks. We give a negative answer. In particular it is not possible to consider peaks which repel each other with an infinitely small speed wrt the period.\\
Our main result is the following uniqueness result, up to a translation in $x_1$, for the small values of $\ep$.\\

\begin{theorem}\label{uniqueDancer}
Let $u$ be a solution of (\ref{eq:u}) that admits the points $\frac{a_{\ep}^1}{\ep}$,...,$\frac{a_{\ep}^k}{\ep}$ in $[-\frac{\pi}{\ep},\frac{\pi}{\ep}[$ as peaks in the sense of the definition \ref{pics}. Then, for $\ep$ small enough there exists $\alpha_{\ep}\rightarrow0$ such that
$$u(x_1,x')=u_D(x_1-\frac{a_{\ep}^1}{\ep}-\alpha_{\ep},x')$$
where $u_D$ is the Dancer solution of period $\frac{2\pi}{k\ep}$.\\
Moreover, if we write $u_D$ as
 \begin{equation}\label{eq:nouveau}
 u_D(x)=\sum_{l\in\Z}U(x_1+{2l\pi\over k\ep},x')+\psi(x)
 \end{equation}
 and if we define 
$$d_x=\hbox{dist}(x,\cup_{l\in\Z}\{({2l\pi\over k\ep},0)\}),$$
then for all $0<\eta<1$ and all $\eta'$ such that $0<\eta'<1$  and $0<\eta'<p-1-\eta$, there exists $C$ independent of $\ep$ such that
 
 \begin{equation}\label{eq:estimationpsi}
 (\vert\psi\vert+\vert\nabla \psi\vert)(x)\leq Ce^{-\eta d_x}e^{-{2\eta'\pi\over k\ep}}({\pi\over k\ep})^{1-N\over2}.
 \end{equation}
\end{theorem}

The most involved part of the proof of Theorem \ref{uniqueDancer} is to prove that the peaks are asymptotically uniformly distributed. More precisely, we will begin with the proof of the following

\begin{proposition}\label{unicite}
Let $u$ be a solution of (\ref{eq:u}) admetting the points $\frac{a_{\ep}^1}{\ep}$,...,$\frac{a_{\ep}^k}{\ep}$ in $[-\frac{\pi}{\ep},\frac{\pi}{\ep}[$ as peaks.\\
Then we have necessarily
 \begin{equation}\label{eq:equirepartis}
 \frac{a_{\ep}^{i+1}-a_{\ep}^i}{\ep}-\frac{2\pi}{k\ep}\rightarrow0,\quad i=0,...,k.
 \end{equation}

\end{proposition}

 Although this property is not useful here, let us recall that we already know that the even ${2\pi\over\ep}$-periodic solution which verifies (\ref{eq:Dancer})
is unique, for $\ep$ small enough. We can prove it as in \cite{singly}, Proposition 3.6, following ideas of \cite{Dancer2}, p. 969. We will use the same kind of proof in the present paper, to complete the proof of Theorem \ref{uniqueDancer} (section 4). Actually, the present paper gives, in particular, a proof of the uniqueness, for $\ep$ small enough, of the ${2\pi\over\ep}$-periodic solution which verifies (\ref{eq:Dancer}) and which admits a maximum at $x=0$, without assuming that this solution is even in $x_1$.\\

 In \cite{Malchiodi}, part 3, Malchiodi gives a construction of a periodic solution with one peak, using a Lyapunov-Schmitt method. 
 
Let us quote the following
\begin{proposition}(Malchiodi,  \cite{Malchiodi}, Corollary 3.2.)\label{Malchiodi}
For $1<p<\frac{N+2}{N-2}$, there exists a solution of (\ref{eq:u}), even in $x_1$, of the form
\begin{equation}\label{eq:formev}
v=\sum_{i\in\Z}U(x_1+i\frac{2\pi}{\ep},x')+\overline w
\end{equation}
where
$$\|\overline w\|_{H^1(]-\frac{\pi}{\ep},\frac{\pi}{\ep}[\times\R^{N-1})}\rightarrow0$$
and
\begin{equation}\label{eq:vu}
\vert\overline w(x)\vert+\vert\nabla\overline w(x)\vert\leq Ce^{-\frac{\pi}{\ep}(1+\xi_0)}e^{-\eta_0\hbox{dist}(x,\cup_{l\in\Z}\{(\frac{2l\pi}{\ep},0)\})}
\end{equation}
for some $\xi_0>0$ and for some $\eta_0>0$.
\end{proposition}

Let us remark that this solution is the Dancer solution, in consideration of the uniqueness of the even ${2\pi\over\ep}$-periodic solution which verifies (\ref{eq:Dancer}). 
In that previous work, the functions are assumed to be even in $x_1$. In ours, we have to overcome some difficulties arising from the lack of evenness. Finally, we prove that the solution is even. We do not use the  Malchiodi's result in the present paper, but it is worth to mention it, for building solutions verifying (\ref{eq:conv}) for the most general configuration of peaks was the initial question of the present work.  \\
 
In the course of the proof of Theorem \ref{unicite}, we will consider an approximate solution of (\ref{eq:u}).

 Let us denote 
 $$ U_{i,l}=U(x_1-{a_{\ep}^i+2\pi l\over\ep},x'),\quad i=1,...,k,\quad l\in\Z,$$
 then
 $$U_{i,0}=U_i.$$
Let us define
$$v_i= \sum_{l\in\Z}U_{i,l}\quad\hbox{and}\quad
\overline u_{\ep}=\sum_{i=1}^kv_i.$$
We will consider the linearized operator about this approximate solution, namely
$$\LL=-\Delta+1-p\overline u^{p-1}_{\ep}.$$
 We will study the vector space associated to the eigenvalues which tend to 0.\\
The operator $(-\Delta+I)^{-1}\LL$ is an operator of $H^1(\frac{S^1}{\ep}\times\R^{N-1})$ into itself of the form $id-{\mathcal K}$, where ${\mathcal K}$ is a compact operator. So $(-\Delta+I)^{-1}\LL$ is a Fredholm operator of index 0.\\
We consider the eigenvalues of the operator $\LL$, in the following sense
\begin{equation}\label{eq:sensvp}
\hbox{there exists $\xi\in H^1(\frac{S^1}{\ep}\times\R^{N-1})$, $\xi\neq0$, verifying
$\LL\xi=\lambda(-\Delta+1)\xi.$}
\end{equation}
The operator $\LL$ has a countably infinite discrete set of eigenvalues, $\lambda_i$, $i=1,2...$. If we designate by $V_i$ the eigenspace corresponding to $\lambda_i$, by $H^1$ the space $H^1(\frac{S^1}{\ep}\times\R^{N-1})$ and by $L^2$ the space $L^2(\frac {S^1}{\ep}\times\R^{N-1})$, then
\begin{equation}\label{eq:valeurspropres}
\lambda_i=\inf\{\frac{<\LL u, u>_{L^2}}{<u,u>_{H^1}},u\neq0,<u,v>_{H^1}=0,\forall v\in V_1\oplus...\oplus V_{i-1}\},\quad\hbox{for $i\geq2$},
\end{equation}
and 
$$\lambda_1=\inf\{\frac{<\LL u, u>_{L^2}}{<u,u>_{H^1}},u\neq0\},$$
(see \cite{GT}).
Let us quote the following result concerning the eigenvalues of the operator $-\Delta+1-pU^{p-1}$ (with the definition above, with $\R^N$ instead of $\frac{S^1}{\ep}\times\R^{N-1}$ and when $k=1$ and $a_{\ep}^1=0$)
\begin{theorem}\label{connu}
The first eigenvalue of $-\Delta+1-pU^{p-1}$ in $\R^N$ is $1-p$. The eigenspace associated with the eigenvalue 0 is spanned by the eigenvectors $\frac{\partial U}{\partial x_j}$, $j=1,...,N$.  
\end{theorem}
This theorem follows from \cite{BL}.\\
Let us define
\begin{equation}\label{eq:defsigmai}
\sigma_i={1\over2}\hbox{dist}(\frac{a_{\ep}^i}{\ep},\cup_{j\neq i,j=0,...,k+1}\{\frac{a_{\ep}^j}{\ep}\})\quad i=1,...,k.
\end{equation}
and
$$\underline\sigma=\min_{i=1}^k\sigma_i.$$

Let us summarize the properties of the eigenfunctions of $\LL$ in the following

\begin{theorem}\label{basedevect}
(i)The eigenvalues of $\LL$ are less than 1. There exists a sequence $(\ep_m)_{m\in\N}$ such that each eigenvalue of $\LL$ tends either to 1 or to an eigenvalue of $-\Delta+1-pU^{p-1}$ as $\ep_m\rightarrow 0$. \\

(ii) Let $F$ be the vector space associated with the eigenvalues tending to 0. Then the dimension of $F$ is $k$ and
$F$ is spanned by  $k$ eigenvectors $\varphi_i$, $i=1,...,k$ such that there exist $k$ real numbers $\alpha_i\neq0$, independent of $\ep$, verifying
\begin{equation}\label{eq:orthonormé}<\varphi_i,\varphi_j>_{H^1(\frac{S^1}{\ep}\times\R^{N-1})}=0\quad i\neq j\quad;\quad\|\varphi_i\|_{{\infty}}=1
\end{equation}
and
\begin{equation}\label{eq:normetheo}\|\varphi_i-\alpha_i \frac{\partial v_i}{\partial x_1}\|_{L^q(\frac{S^1}{\ep}\times\R^{N-1})}+\|\nabla( \varphi_i-\alpha_i \frac{\partial v_i}{\partial x_1})\|_{L^q(\frac{S^1}{\ep}\times\R^{N-1})}\rightarrow0
\end{equation}
for all $1\leq q\leq\infty$.\\
\end{theorem}

 The paper is organized as follows. In section 2, we study the eigenvectors associated with the eigenvalues of the operator $\LL$ which tend to 0 and we give the proof of Theorem \ref{basedevect}.  In section 3, we use a Lyapunov-Schmitt method to give the proof of Proposition \ref{unicite}. In section 4, we conclude the proof of Theorem \ref{uniqueDancer}.
 
 In sections 2, 3 and 4 we will refer to some technical results, which are reported in the appendix (section 5).

\section{An analysis of the eigenvalues.}

In this part, we prove the theorem \ref{basedevect}. \\

Proof of (i).\\

 Let $\varphi$  be such that
$$\LL \varphi=\lambda(-\Delta \varphi+\varphi)\quad\hbox{in ${S^1\over\ep}\times\R^{N-1}$}.$$
We suppose that there exists $c$ such that $\varphi(c)\rightarrow 1$ and that $\|\varphi\|_{\infty}=1$.
 We denote 
$$\phi(x)=\varphi(x+c).$$
By standard elliptic estimates, there exists a subsequence such that $\phi\rightarrow\overline \phi$ uniformly on compact sets of $\R^N$.
Let us suppose that $\lambda\not\rightarrow 1$. First, if $\vert c-(\frac{a_{\ep}^i}{\ep},0)\vert\rightarrow+\infty$ for all $i$, then 
$$-\Delta\overline \phi+\overline \phi=0\quad\hbox{ in $\R^N$}\quad;\quad\|\overline \phi\|_{\infty}=\overline \phi(0)=1.$$
This is in contradiction with the maximum principle, so this case does not occur.
So there exists $\overline c$ and $i$ such that $(c-(\frac{a_{\ep}^i}{\ep},0))\rightarrow\overline c$. Then, 
\begin{equation}\label{eq:U}
(-\Delta+1-pU^{p-1}(x+\overline c))\overline \phi=\overline\lambda(-\Delta\overline \phi+\overline \phi)\quad\hbox{ in $\R^N$}
\end{equation}
and $\overline \phi$ is non zero and even in $x'$. \\
Let $x\in\R^N$ be given. Since $c_1-\frac{a^i_{\ep}}{\ep}\rightarrow\overline c_1$, then we have, for $\ep$ small enough
$$-\tilde\sigma_{i-1}<x_1+c_1-\frac{a^i_{\ep}}{\ep}<\tilde\sigma_i,$$
where
$$\tilde \sigma_j=\frac{a_{\ep}^{i+1}-a_{\ep}^{i}}{2\ep}\quad\hbox{for $j=0,...,k$}.$$
In view of (\ref{eq:fp}), we deduce that
$$\vert\phi(x)\vert\leq Ce^{-\eta\vert x+c-(\frac{a^i_{\ep}}{\ep},0)\vert}
$$
where $C$ is independent of $\ep$. Letting $\ep\rightarrow0$, we obtain
$$\vert\overline\phi(x)\vert\leq Ce^{-\eta\vert x+\overline c\vert}
$$
and this is true for any $x\in\R^N$. Thus $\overline\lambda$ is an eigenvalue of $-\Delta+1-pU^{p-1}$ in the sense of (\ref {eq:sensvp}). \\
  By a diagonal process, we can construct a subsequence $(\ep_m)$ such that any eigenvalue of $\LL$ which does not tend to 1 converges to an eigenvalue of $-\Delta+1-pU^{p-1}$.\\

Proof of (ii).\\

 We divide the proof into three parts.\\
Firstly, let us prove that if $\varphi\in F\backslash\{0\}$, $\|\varphi\|_{\infty}=1$, then there exists $I\subset\{1,...,k\}$ and some real numbers $\beta_i\neq0$ and independent of $\ep$ such that
\begin{equation}\label{eq:norme}
\|\varphi -\sum_{i\in I}\beta_i{\partial v_i\over\partial x_1} \|_{\infty}+\|\nabla(\varphi -\sum_{i\in I}\beta_i{\partial v_i\over\partial x_1})\|_{\infty}\rightarrow0.
\end{equation}
We follow the proof of (i) to get (\ref{eq:U}) with $\overline\lambda=0$. 
There exists some real number $\beta\neq0$ such that
$$\overline \phi(x)=\beta\frac{ \partial U}{\partial x_1}(x+\overline c).$$
 We get 
$$\varphi(x+c)-\beta\frac{\partial U}{\partial x_1}(x+\overline c)\rightarrow0\quad\hbox{uniformly on compact sets},$$
that is
$$\varphi(x+c)-\beta\frac{\partial U}{\partial x_1}(x+c-(\frac{a_{\ep}^i}{\ep},0))\rightarrow0\quad\hbox{uniformly on compact sets},$$
that leads to

$$\quad\varphi(x)-\beta\frac{\partial U_i}{\partial x_1}(x)\rightarrow0\quad\hbox{ uniformly for $x$ such that $(x_1-\frac{a_{\ep}^i}{\ep})$ is bounded}.$$
Finally for each $i$, either there exists $\alpha_i\neq0$ such that
$$(\varphi-\beta_i\frac{\partial v_i}{\partial x_1})\rightarrow 0\quad\hbox{ uniformly for $x$ such that $(x_1-\frac{a_{\ep}^i}{\ep})$ is bounded}$$
or
 $$\varphi\rightarrow 0\quad\hbox{ uniformly for $x$ such that $(x_1-\frac{a_{\ep}^i}{\ep})$ is bounded}.$$ Moreover, the first case occurs for at least one $i$. 
 We can use Proposition \ref{barrierefonctionpropre} to conclude that there exists $I\subset\{1,...,k\}$ and $\beta_i\neq 0$ and independent of $\ep$ such that
\begin{equation}\label{eq:I1}\|\varphi-\sum_{i\in I}\beta_i\frac{\partial v_i}{\partial x_1}\|_{\infty}\rightarrow0.
\end{equation}

We deduce that
\begin{equation}\label{eq:I2}
\|\nabla(\varphi-\sum_{i\in I}\beta_i\frac{\partial v_i}{\partial x_1})\|_{\infty}\rightarrow0
\end{equation}

 by standard elliptic arguments.\\


Secondly, let us assume that $F\neq\{0\}$.
The subsequence $(\ep_m)$ being defined in (i),
we define a finite set $ J$ and eigenvalues $\lambda_j$, $j\in J$ such that $\lambda_j(\ep_m)\rightarrow0$. Let $\varphi_j$ be an eigenvector associated with $\lambda_j$. 
Let us assume that $(\varphi_i)_{i\in  J}$ is a basis of $F$ which verifies
 $$<\varphi_i,\varphi_j>_{H^1([-{\pi\over\ep},{\pi\over\ep}]\times\R^{N-1})}=0\quad i\neq j\quad;\quad\|\varphi_i\|_{L^{\infty}([-{\pi\over\ep},{\pi\over\ep}]\times\R^{N-1})}=1.$$
 
We write
\begin{equation}\label{eq:theo}
\frac{\partial v_{i}}{\partial x_1}=
\sum_{l\in J}c^i_l\varphi_l+\xi_i\quad;\quad
<\xi_i,\varphi_l>_{H^1}=0\hbox{ for  $l\in J$}.
\end{equation}
We have
\begin{equation}\label{eq:Lxi}
\LL\xi_i=\LL\frac{\partial v_i}{\partial x_1}
-\sum_lc_l^i\lambda_l(-\Delta\varphi_l+\varphi_l).
\end{equation}
In view of (\ref{eq:xiinfini}), we have
$$\|\xi_i\|_{\infty}\leq C\|\LL\xi_i\|_{\infty}.$$
Since
$$\|\LL\frac{\partial v_i}{\partial x_1}\|_{\infty}\leq p\sum_{l\in\Z}\|(\overline u_{\ep}^{p-1}-U_{i,l}^{p-1})\frac{\partial U_{i,l}}{\partial x_1}\|_{\infty}\rightarrow0,$$
we deduce that
\begin{equation}\label{eq:xitend0}
\|\xi_i\|_{\infty}\rightarrow0.
\end{equation}
 Now, since $\|\frac{\partial v_i}{\partial x_1}\|_{\infty}\not\rightarrow0$, we have that $F\neq\{0\}$. 
 
Thirdly, since $c^i_l=<\frac{\partial v_{i}}{\partial x_1},\varphi_l>_{H^1}$, we may define
$${d^i_l}=\lim_{\ep\rightarrow0}c^i_l.$$
We remark that for all $i$, there exists $l$ such that ${d^i_l}\neq0.$\\
For $i\neq j$, taking the scalar product in $H^1$ we obtain
 
 $$<\frac{\partial v_{i}}{\partial x_1},\frac{\partial v_{j}}{\partial x_1}>_{H^1}=\sum_{l\in J}c^i_lc^j_l\|\varphi_l\|^2_{H^1}+<\xi_i,\xi_j>_{H^1}.$$
 
In view of Proposition \ref{barrierefonctionpropre} and of (\ref{eq:xitend0}), the Lebesgue Theorem leads to
$$\sum_{l\in J}d^i_ld^j_l=0.$$

So, we obtain an orthogonal family $(\tilde\varphi_1,...,\tilde\varphi_k)$ of $F$ defined by
$$\tilde\varphi_i={\alpha_i}{\sum_{l\in J}d_l^i\varphi_l},\quad{i=1,...,k}$$
where
$$(\alpha_i)^{-1}=\|\sum_{l=1}^kd_l^i\varphi_l\|_{\infty}.$$
Then we have
\begin{equation}\label{eq:tildeph}
\|\tilde\varphi_i-\alpha_i\frac{\partial v_i}{\partial x_1}\|_{\infty}+\|\nabla(\tilde\varphi_i-\alpha_i \frac{\partial v_i}{\partial x_1})\|_{\infty}\rightarrow0\quad\hbox{$i=1,...,k$}.
\end{equation}
Now, let us prove that $F$ is spanned by $(\tilde\varphi_1,...,\tilde\varphi_k)$.\\
 Otherwise, let $\varphi\in F$, $<\varphi,\tilde\varphi_i>=0$, $i=1,...,k$,
$\|\varphi\|_{\infty}=1$. Using (\ref{eq:I1}) and (\ref{eq:tildeph}), we obtain
$$\|\varphi-\sum_{i\in I}\frac{\beta_i}{\alpha_i}\tilde\varphi_i\|_{\infty}+\|\nabla(\varphi-\sum_{i\in I}\frac{\beta_i}{\alpha_i}\tilde\varphi_i)\|_{\infty}\rightarrow0.$$
By Proposition \ref{barrierefonctionpropre}, we deduce that
$$\|\varphi-\sum_{i\in I}\frac{\beta_i}{\alpha_i}\tilde\varphi_i\|_{H^1}\rightarrow0$$
and we are led to
 $$\|\varphi\|^2_{H^1}+\sum_{i\in I}(\frac{\beta_i}{\alpha_i})^2\|\tilde\varphi_i\|^2_{H^1}\rightarrow0,$$ that  contradicts (\ref{eq:normes}).\\
Consequently, $(\tilde\varphi_1,...,\tilde\varphi_k)$ is an orthogonal basis of $F$. Denoting again this new basis by $(\varphi_1,...,\varphi_k)$, we get a basis of $F$ verifying (\ref{eq:orthonormé}). 
Moreover, by Proposition \ref{barrierefonctionpropre}, $\varphi_i$ is bounded in $L^q$ and $\|\xi\|_{L^q}\rightarrow0$ for all $q\geq1$, that gives (\ref{eq:normetheo}).\\

\section{ The Lyapunov-Schmitt reduction.}

In this part, we prove the proposition \ref{unicite}.\\

Let us define
$${\mathcal M}(u)=-\Delta u+u-u_+^p.$$

Let $k$ real numbers $\delta_1$,...,$\delta_k$ and $v\in H^1(\frac{S^1}{\ep}\times\R^{N-1})$ be given. Let us suppose 
that
    \begin{equation}\label{eq:solution}
  u=\overline u_{\ep}+v+\sum_{i=1}^k \delta_i\varphi_i;\quad <v,\varphi_i>_{H^1}=0,\quad i=1,...,k
\end{equation}          
is a solution of (\ref{eq:u}) and that 
  $$\| v\|_{\infty}+\sum_{i=1}^k\vert \delta_i\vert\rightarrow0.$$  
 
 We define
\begin{equation}\label{eq:h}
h=-\mathcal{M}(\overline u_{\ep}+v+\sum_{i=1}^k \delta_i\varphi_i)+\LL (v+\sum_{i=1}^k \delta_i\varphi_i).
\end{equation}
Then $v$ and $\delta_1$,...,$\delta_k$ are such that
$$\LL(v+\sum_{i=1}^k \delta_i\varphi_i)=h.$$
 We denote
 $$h=h^{\bot}+h^{ \top}\quad,\quad h^{ \top}\in\hbox{Vect}\{\varphi_1,...,\varphi_k\}\quad,\quad h^{ \bot}\in (\hbox{Vect}\{\varphi_1,...,\varphi_k\})^{\bot},$$ 
 relatively to the Hilbert space $H^1([-\frac{\pi}{\ep},\frac{\pi}{\ep}]\times\R^{N-1})$. First, $v$ is a ${2\pi\over\ep}$-periodic solution of the equation
\begin{equation}\label{eq:pointfixe}
\left\{\begin{array}{rl}
&\LL v=h^{\bot}\\
&<v,\varphi_i>_{H^1}=0,\quad i=1,...,k.
\end{array}  
\right.                                       
\end{equation}
Then $(\delta_1,...,\delta_k)$ verifies
$$\LL(\sum_{i=1}^k \delta_i\varphi_i)=h^{\top}.$$
To begin with, we have
\begin{lemma}\label{Mappr}
\begin{equation}\label{eq:Mappr}\|{\mathcal M}(\overline u_{\ep})
\|_{\infty}\leq C e^{-2\underline\sigma}\underline\sigma^{1-N\over2}.
\end{equation}
\end{lemma}

{\bf Proof.}
\begin{equation}\label{eq:Mest}
{\mathcal M}(\overline u_{\ep})=\sum_{i,l} U_{i,l}^p-(\sum_{i,l} U_{i,l})^p.
\end{equation}
We define, for $i=0,...,k+1$
\begin{equation}\label{eq:defomega}
\Omega_{i}=\{x\in[-{\pi\over\ep},{\pi\over\ep}]\times\R^{N-1};\hbox{dist}( x,\cup_{l=0}^{k+1}\{(\frac{a_{\ep}^l}{\ep},0)\}=\vert x-(\frac{a_{\ep}^i}{\ep},0)\vert\}.
\end{equation}
We have
$$[-{\pi\over\ep},{\pi\over\ep}]\times\R^{N-1}=\bigcup_{i=0}^{k+1}\Omega_i.$$

For all $x\in\Omega_i$, for all $j\neq i$, we have together
$$\sum_{l\in\Z}U_{j,l}+\sum_{l\in\Z^{\star}}U_{i,l}\leq Ce^{-2\sigma_i}\sigma_i^{1-N\over2} e^{\vert x-\frac{a_{\ep}^i}{\ep}\vert}
$$
and
$$\sum_{l\in\Z}U_{j,l}+\sum_{l\in\Z^{\star}}U_{i,l}\leq Ce^{-\sigma_i}\sigma_i^{1-N\over2} .
$$
Then, by Lemma \ref{taylor} we write in $\Omega_i$
$${\mathcal M}(\overline u_{\ep})=-p U_i^{p-1}(\sum_{j\neq i;l\in\Z}U_{j,l}+\sum_{l\in\Z^{\star}}U_{i,l})+O(\sum_{j\neq i;l\in\Z}U_{j,l}+\sum_{l\in\Z^{\star}}U_{i,l})^2+\sum_{j\neq i;\l\in\Z}U^p_{j,l}$$
while $p-1\geq1$. We easily deduce the proof of the Lemma.\\

  We have the following 
\begin{proposition}\label{normev}
 Let $v$ be a solution of (\ref{eq:pointfixe}). Then there exists  $C$ independent of $\ep$ such that
\begin{equation}\label{eq:normev}
\hbox{if $p>2$}\quad
\| v\|_{H^1}\leq C(e^{-2\underline\sigma}\underline\sigma^{1-N\over2} +\sum_{i=1}^k \vert\delta_i\vert^{2}) 
\end{equation}
$$\hbox{ if $p=2$ $\forall\eta\in]0,1[$},\quad \| v\|_{H^1}\leq C(e^{-2\eta\underline\sigma} +\sum_{i=1}^k \vert\delta_i\vert^{2}) $$
and for all $p$
\begin{equation}\label{eq:normevinfini}
\|v\|_{\infty}+\|\nabla v\|_{\infty}
\leq C( e^{-2\underline\sigma}\underline\sigma^{1-N\over2} +\sum_{i=1}^k \vert\delta_i\vert^{2}).
\end{equation}

\end{proposition}

{\bf Proof.}  

We write
$$h=\Delta\overline u_{\ep}-\overline u_{\ep}+(\overline u_{\ep}+v+\sum\delta_i\varphi_i)^p_+-p\overline u^{p-1}_{\ep}(v+\sum\delta_i\varphi_i)$$
and Lemma \ref{taylor}
gives
$$\vert(\overline u_{\ep}+v+\sum\delta_i\varphi_i)^p_+-\overline u^p_{\ep}-p\overline u^{p-1}_{\ep}(v+\sum\delta_i\varphi_i)\vert\leq C\vert v+\sum\delta_i\varphi_i\vert^{2}.$$
We deduce that
\begin{equation}\label{eq:estih}
\vert h+{\mathcal M}(\overline u_{\ep})\vert\leq C\vert v+\sum\delta_i\varphi_i\vert^{2}.
\end{equation}

Since
$$h^{\bot}=h-\sum_{i=1}^k\frac{<h,\varphi_i>_{L^2}}{\|\varphi_i\|^2_{H^1}}(-\Delta\varphi_i+\varphi_i),$$
then
$$\|h^{\bot}\|_{L^2}\leq C\|h\|_{L^2}
\quad\hbox{and}\quad\|h^{\bot}\|_{\infty}\leq C\|h\|_{\infty}.
$$
Now, we use (\ref{eq:xiinfini}) to obtain
$$\|v\|_{\infty}\leq C\|h^{\bot}\|_{\infty}.
$$
Using (\ref{eq:Mappr}), we deduce the estimate
\begin{equation}\label{eq:vtransit}
\|v\|_{\infty}\leq C( e^{-2\underline\sigma}\underline\sigma^{1-N\over2} +\sum_{i=1}^k \vert\delta_i\vert^{2})
\end{equation}
and the estimate (\ref{eq:normevinfini}) follows in the standard way.\\
We have also by (\ref{eq:minorant})
$$\|v\|_{H^1}\leq C\|h^{\bot}\|_{L^2}.$$
We write

$$\|h\|^2_{L^2}\leq C(\|{\mathcal M}(\overline u_{\ep})\|^2_{L^2}+\|v\|_{L^4}^4+\sum_j\vert\delta_j\vert^4).
$$

Using (\ref{eq:vtransit}), we deduce, for $\ep$ small enough
$$\|v\|_{H^1}\leq C(\|{\mathcal M}(\overline u_{\ep})\|_{L^2}+\sum_j\vert\delta_j\vert^2).$$

Now we use (\ref{eq:Mest}) and Lemma \ref{taylor} to obtain
$$\|{\mathcal M}(\overline u_{\ep})\|^2_{L^2}\leq\sum_{i=0}^{k+1}\int_{\Omega_i} pU_i^{2(p-1)}(\sum_{j\neq i;l\in\Z}U_{j,l}+\sum_{l\in\Z^{\star}}U_{i,l})^2+ C\int_{\Omega_i}(\sum_{j\neq i;l\in\Z}U_{j,l}+\sum_{l\in\Z^{\star}}U_{i,l})^4.
$$
 In order to estimate the first integral, we use Proposition \ref{convolution2} with $a=2(p-1)$ and $b=2$, as long as $p>2$, and with $a=2$ and $b=2\eta$, $0<\eta<1$, if $p=2$. Finally, we obtain
\begin{equation}\label{eq:M2}
\|{\mathcal M}(\overline u_{\ep})\|_{L^2}\leq C(e^{-2\underline\sigma}\underline\sigma^{1-N\over2})\quad\hbox{if $p>2$}
\end{equation}
and for all $\eta\in]0,1[$
\begin{equation}\label{eq:M2p2}
\|{\mathcal M}(\overline u_{\ep})\|_{L^2}\leq C(e^{-2\eta\underline\sigma}\underline\sigma^{1-N\over2})\quad\hbox{if $p=2$}.
\end{equation}
We deduce (\ref{eq:normev}).\\
We have proved the proposition.\\
\begin{proposition}\label{deltai}
Let $u$ be given as in Proposition \ref{unicite} and let $\delta_1$,....,$\delta_k$ be defined in (\ref{eq:solution}). We can possibly replace the $k$ given points $\frac{a_{\ep}^1}{\ep}$,....,$\frac{a_{\ep}^k}{\ep}$ by $k$ points $\frac{b_{\ep}^1}{\ep}$,....,$\frac{b_{\ep}^k}{\ep}$
verifying
$$ \frac{a_{\ep}^i}{\ep}-\frac{b_{\ep}^i}{\ep}\rightarrow0\quad\hbox{as $\ep\rightarrow0$}$$
in order to have
$$\delta_i=0,\quad i=1,...,k.$$
\end{proposition}

{\bf Proof.} Let 
$$u=\overline u_{\ep}+\sum_{j=1}^k\delta_j\varphi_j+v$$
be the given solution of (\ref{eq:u}).\\
Let us give $(\alpha_1,...,\alpha_k)$ depending on $\ep$, such that  $(\alpha_1,...,\alpha_k)\rightarrow0$. We can replace the points $ \frac{a_{\ep}^i}{\ep}$ by the points $ \frac{a_{\ep}^i}{\ep}+\alpha_i$. In other words, we write
$$ u=\tilde u_{\ep}+\sum_{j=1}^k\tilde\delta_j\tilde\varphi_j+\tilde v$$
where 
$$\tilde u_{\ep}(x)=\sum_{j=1}^k\sum_{l\in\Z}U_{j,l}(x_1-\alpha_j,x')$$
$$\hbox{and}\quad\tilde\varphi_j, \quad j=1,...,k$$ are the eigenfunctions corresponding to the eigenvalues tending to 0, for the configuration of points $\frac{a_{\ep}^j}{\ep}+\alpha_j$, and $$<\tilde v,\tilde\varphi_j>=0, \quad j=1,...,k.$$
Substracting the expressions of $u$ and performing the scalar product in $H^1$ by $\varphi_i$, we get
\begin{equation}\label{eq:tildedelta}
\tilde\delta_i\|\varphi_i\|^2_{H^1}+\sum_j\tilde\delta_j<\tilde\varphi_j-\varphi_j,\varphi_i>_{H^1}=<\overline u_{\ep}-\tilde u_{\ep},\varphi_i>_{H^1}+\delta_i\|\varphi_i\|^2_{H^1}+<v-\tilde v,\varphi_i>_{H^1}.
\end{equation}
First we remark that
 we have $$<v-\tilde v,\varphi_i>=<\tilde v,\tilde\varphi_i-\varphi_i>$$
while by (\ref{eq:normev})
$$\|\tilde v\|_{H^1}\leq C(e^{-2\eta\underline\sigma}\underline\sigma^{1-N\over2}+\sum_j\tilde\delta_j^2),$$
with $\eta=1$, for $p>2$,
thus
\begin{equation}\label{eq:v-tilde}
\vert<v-\tilde v,\varphi_i>\vert\leq C(e^{-2\eta\underline\sigma}\underline\sigma^{1-N\over2}+\sum_j\tilde\delta_j^2).
\end{equation}
Moreover
$$\vert<\overline u_{\ep}-\tilde u_{\ep},\varphi_i>\vert\leq C\sum_j\vert\alpha_j\vert$$

and, as a consequence of (\ref{eq:normetheo})
$$\|\tilde\varphi_j-\varphi_j\|_{H^1}\rightarrow0.$$
Thanks to (\ref{eq:tildedelta}), we deduce that 
\begin{equation}\label{eq:pre}
\sum_i\vert\tilde\delta_i\vert\leq C(\sum_j\vert\delta_j\vert+\sum_j\vert\alpha_j\vert+e^{-2\eta\underline\sigma}\underline\sigma^{1-N\over2}).
\end{equation}
and consequently

$$\|\tilde v\|_{H^1}\leq C(e^{-2\eta\underline\sigma}\underline\sigma^{1-N\over2}+\sum_j\vert\delta_j\vert^2+\sum_j\vert\alpha_j\vert^2),
$$
thus
\begin{equation}\label{eq:pre}
\vert <v-\tilde v,\varphi_i>\vert\leq C(e^{-2\eta\underline\sigma}\underline\sigma^{1-N\over2}+\sum_j\vert\delta_j\vert^2+\sum_j\vert\alpha_j\vert^2)
\end{equation}
for some $C$ independent of $\ep$
.\\
Now let us prove that we can choose $(\alpha_1,...,\alpha_k)$ such that
$$<\overline u_{\ep}-\tilde u_{\ep},\varphi_i>+\delta_i\|\varphi_i\|^2+<v-\tilde v,\varphi_i>=0.$$

We define
$$\mathcal{F}(\alpha_1,...,\alpha_k)=(<\overline u_{\ep}-\tilde u_{\ep},\varphi_i>)_{i=1,...,k}.$$
This definition gives, for $i$ and $j=1,...,k$
$${\partial \mathcal{F}_i\over\partial \alpha_j}=\sum_{l\in\Z}\int_{[-\frac{\pi}{\ep},\frac{\pi}{\ep}]\times\R^{N-1}}\frac{\partial U_{j,l}}{\partial x_1}(x_1-\alpha_j,x')\varphi_i(x)dx+\sum_{l\in\Z}\int_{[-\frac{\pi}{\ep},\frac{\pi}{\ep}]\times\R^{N-1}} \nabla\frac{ \partial U_{j,l}}{\partial x_1}(x_1-\alpha_j,x').\nabla \varphi_i(x)dx
.$$
We deduce that, as $\ep\rightarrow0$
$${\partial \mathcal{F}_i\over\partial \alpha_i}(0)\rightarrow \|\frac{\partial U}{\partial x_1}\|^2_{H^1}
\quad\hbox{
and }\quad
{\partial \mathcal{F}_i\over\partial \alpha_j}(0)\rightarrow 0\quad\hbox{ for $j\neq i$}.$$
Thus $d\mathcal{F}(0)$ is an isomorphism, for $\ep$ small enough.\\
Let us define $\alpha=(\alpha_1,...,\alpha_k)$.
We have to solve
$$\mathcal{F}(\alpha)+(\delta_i\|\varphi_i\|^2+<v-\tilde v,\varphi_i>)_{i=1,...,k}=0.$$
 We define
$$Q(\alpha)=\mathcal{F}(\alpha)-\mathcal{F}(0)-d\mathcal{F}(0)(\alpha)$$
and
$$\mathcal{G}(\alpha)=(-d\mathcal{F}(0))^{-1}(Q(\alpha)+(\delta_i\|\varphi_i\|^2+<v-\tilde v,\varphi_i>)_{i,...,k}).$$
Since we have together
$$\vert Q(\alpha)\vert=O(\vert\alpha\vert^2)$$
and (\ref{eq:pre}),
we can use the Brouwer fixed point Theorem in a standard way. We find a real number $R$,
$$R\leq C(e^{-2\eta\underline\sigma}\underline\sigma^{1-N\over2}+\sum_j\vert\delta_j\vert)$$
 such that 
$$(\vert\alpha\vert\leq R)\Rightarrow\vert\mathcal{G}(\alpha)\vert\leq R).$$
So we find $\alpha$, $\vert\alpha\vert\leq R$, such that $\mathcal{G}(\alpha)=\alpha$, that is
$$<\overline u_{\ep}-\tilde u_{\ep},\varphi_i>+\delta_i\|\varphi_i\|^2+<v-\tilde v,\varphi_i>=0.$$
Returning to (\ref{eq:tildedelta}), we deduce that $\tilde\delta_i=0$, $i=1,...,k$.\\

Now, we suppose that $\delta_i=0$, $i=1,...,k$.\\
Let $d_i$ be defined by
$$h^{\top}=\sum_{i=1}^kd_i(-\Delta\varphi_i+\varphi_i).$$
On one hand, since 
$$h^{\top}=\LL(\sum_{i=1}^k\delta_i\varphi_i)=0$$
we have
\begin{equation}\label{eq:deq0}
d_i=0\quad\hbox{for all $i$.}
\end{equation}
On the other hand, we have the following

\begin{proposition}\label{di}

For $i=1,...,k$
\begin{equation}\label{eq:di}
d_i={p\over\|\varphi_i\|^2_{H^1}}\int_{\Omega_i}U_i^{p-1}\sum_{j\neq i}v_j{\partial U_i\over\partial x_1}dx
+o(e^{-2\underline\sigma}\underline\sigma^{1-N\over2}).
\end{equation}
\end{proposition}

{\bf Proof.} We suppose that $k\geq 2$, otherwise the proposition \ref{unicite} is irrelevant. We have
$$d_i={1\over\|\varphi_i\|^2_{H^1}}\int_{\frac{S^1}{\ep}\times\R^{N-1}}(h\varphi_i)dx.$$

\begin{equation}\label{eq:di2}
d_i=(\int_{\frac{S^1}{\ep}\times\R^{N-1}}(h+\mathcal{M}(\overline u_{\ep}))\varphi_idx-\int_{\frac{S^1}{\ep}\times\R^{N-1}}\mathcal{M}(\overline u_{\ep})\varphi_i dx){1\over\|\varphi_i\|^2_{H^1}}
\end{equation}
The coefficient $\|\varphi_i\|_{H^1}$ does not matter, thanks to (\ref{eq:normes}). Since $\delta_i=0$, $i=1,...,k$, we deduce from (\ref{eq:normevinfini}) and (\ref{eq:estih}) that
\begin{equation}\label{eq:1er}
\int_{\frac{S^1}{\ep}\times\R^{N-1}}(h+\mathcal{M}(\overline u_{\ep}))\varphi_idx
=O(e^{-4\eta\underline\sigma}\underline\sigma^{1-N}).
\end{equation}
Now let us estimate the second integral of (\ref{eq:di2}), for $i=1,...,k$.\\ 
Without loss of generality, we let $i=1$.  We write
$$\int_{]-\frac{\pi}{\ep},\frac{\pi}{\ep}[\times\R^{N-1}}\mathcal{M}(\overline u_{\ep})\varphi_1dx=\int_{]-\frac{\pi}{\ep},\frac{\pi}{\ep}[\times\R^{N-1}}\mathcal{M}(\overline u_{\ep})(\varphi_1-\alpha_1\frac{\partial v_1}{\partial x_1})dx+\alpha_1\int_{]-\frac{\pi}{\ep},\frac{\pi}{\ep}[\times\R^{N-1}}\mathcal{M}(\overline u_{\ep})\frac{\partial v_1}{\partial x_1}dx.$$
The estimate (\ref{eq:Mappr}) gives directly

\begin{equation}\label{eq:reste1}
\int_{]-\frac{\pi}{\ep},\frac{\pi}{\ep}[\times\R^{N-1}}\mathcal{M}(\overline u_{\ep})(\varphi_1-\alpha_1\frac{\partial v_1}{\partial x_1})dx=o(e^{-2\underline\sigma}\underline\sigma^{1-N\over2}),
\end{equation}
since 
$$\|\varphi_1-\alpha_1\frac{\partial v_1}{\partial x_1}\|_{L^1}\rightarrow0.$$
Now we write
$$\int_{]-\frac{\pi}{\ep},\frac{\pi}{\ep}[\times\R^{N-1}}\mathcal{M}(\overline u_{\ep})\frac{\partial v_1}{\partial x_1}dx=
\int_{\Omega_1}\mathcal{M}(\overline u_{\ep})\frac{\partial v_1}{\partial x_1}dx
+{\sum_{j=0}^{k+1}}_{j\neq 1}p\int_{\Omega_j}\mathcal{M}(\overline u_{\ep})\frac{\partial v_1}{\partial x_1} dx$$

We have 
$$\mathcal{M}(\overline u_{\ep})=\sum_{i,l}U_{i,l}^p-(\sum_{i,l}U_{i,l})^p.$$

By Lemma \ref{taylor} we have, in $\Omega_j$,

$$\mathcal{M}(\overline u_{\ep})=-pU_j^{p-1}(\sum_{i\neq j,l\in\Z}U_{i,l}+\sum_{l\neq0}U_{j,l})+O(\sum_{i\neq j,l\in\Z}U_{i,l}+\sum_{l\neq0}U_{j,l})^2.
$$
Since
$$\vert x-\frac{a_{\ep}^i}{\ep}\vert\geq\sigma_j\quad\hbox{and}\quad\vert x-\frac{a_{\ep}^j+2\pi l}{\ep}\vert\geq\sigma_j\quad\hbox{$\forall i\neq j$ $\forall l\neq0$}$$
then
$$\mathcal{M}(\overline u_{\ep})=-pU_j^{p-1}(\sum_{i\neq j,l\in\Z}U_{i,l}+\sum_{l\neq0}U_{j,l})+O(e^{-2\sigma_j}\sigma_j^{1-N}))\quad\hbox{in $\Omega_j$}.
$$
We get, for $j\neq1$ 
$$
\int_{\Omega_j}\mathcal{M}(\overline u_{\ep})\frac{\partial v_1}{\partial x_1} dx=O(\|\frac{\partial v_1}{\partial x_1}\|^{1\over2}_{L^{\infty}(\Omega_j)})\int_{\Omega_j}\vert\frac{\partial v_1}{\partial x_1}\vert^{1\over2}U_j^{p-1}(\sum_{i\neq j,l}U_{i,l}+\sum_{l\neq0}U_{j,l})dx
+o(e^{-2\sigma_j}\sigma_j^{1-N\over2}).
$$
Now we use Proposition \ref{convolution2}, with $a=p-1+\frac{1}{2}>1$, $b=1$ and with 
$\vert y_0\vert\geq 2\sigma_j$, thanks to
$$\vert\frac{\partial v_1}{\partial x_1}\vert^{1\over2}U_j^{p-1}\leq Ce^{-(p-{1\over2})\vert x-\frac{a_{\ep}^j}{\ep}\vert}\quad\hbox{in $\Omega_j$}$$

 We are led to
 \begin{equation}\label{eq:jnon1}
 \int_{\Omega_j}\mathcal{M}(\overline u_{\ep})\frac{\partial v_1}{\partial x_1} dx=o(e^{-2\sigma_j}\sigma_j^{1-N\over2})\quad\hbox{for $j\neq1$}.
 \end{equation}

Now we write
$$
\int_{\Omega_1}\mathcal{M}(\overline u_{\ep})\frac{\partial v_1}{\partial x_1}  dx=-p\int_{\Omega_1}U_1^{p-1}(\sum_{j\neq1}v_j+\sum_{l\neq0}U_{1,l}){\partial v_1\over\partial x_1}dx
$$
$$+\int_{\Omega_1}{\partial v_1\over\partial x_1}O(\sum_{i\neq 1,l\in\Z}U_{i,l}+\sum_{l\neq0}U_{1,l})^2dx.$$
For any $0<\eta<1$, we have
$$\int_{\Omega_1}{\partial v_1\over\partial x_1}(\sum_{i\neq 1,l\in\Z}U_{i,l}+\sum_{l\neq0}U_{1,l})^2dx=O(\|\sum_{i\neq 1,l\in\Z}U_{i,l}+\sum_{l\neq0}U_{1,l}\|^{1+\eta}_{L^{\infty}(\Omega_1)})\int_{\Omega_1}{\partial v_1\over\partial x_1}(\sum_{i\neq 1,l\in\Z}U_{i,l}+\sum_{l\neq0}U_{1,l})^{1-\eta}dx.$$
We can use Proposition \ref{convolution2} again, with $a=1$, $b=1-\eta$ and $\vert y_0\vert\geq\sigma_1$, to obtain
$$\int_{\Omega_1}{\partial v_1\over\partial x_1}(\sum_{i\neq 1,l\in\Z}U_{i,l}+\sum_{l\neq0}U_{1,l})^{1-\eta}dx=O(e^{-(1-\eta)\sigma_1}\sigma_1^{(1-\eta){1-N\over2}})$$
and then we get
$$O(\|\sum_{i\neq 1,l\in\Z}U_{i,l}+\sum_{l\neq0}U_{1,l}\|^{1+\eta}_{L^{\infty}(\Omega_1)})\int_{\Omega_1}{\partial v_1\over\partial x_1}(\sum_{i\neq 1,l\in\Z}U_{i,l}+\sum_{l\neq0}U_{1,l})^{1-\eta}dx=O(e^{-2\sigma_1}\sigma_1^{1-N})=o(e^{-2\sigma_1}\sigma_1^{1-N\over2})$$
and consequently
\begin{equation}\label{eq:un}
\int_{\Omega_1}\mathcal{M}(\overline u_{\ep})\frac{\partial v_1}{\partial x_1}  dx=-p\int_{\Omega_1}U_1^{p-1}(\sum_{j\neq1}v_j+\sum_{l\neq0}U_{1,l}){\partial v_1\over\partial x_1}dx+o(e^{-2\sigma_1}\sigma_1^{1-N\over2}).
\end{equation}
Moreover, since $k\geq 2$, we have
$$\frac{\pi}{\ep}-\sigma_1\rightarrow+\infty$$
and consequently
$$e^{-2\pi\over\ep}=o(e^{-2\sigma_1}).$$
So, we deduce from Proposition \ref{convolution2} that
$$\int_{\Omega_1}U_1^{p-1}\sum_{l\neq0}U_{1,l}\sum_{l}{\partial U_{1,l}\over\partial x_1}dx=o(e^{-2\sigma_1}\sigma_1^{1-N\over2})$$
and
$$\int_{\Omega_1}U_1^{p-1}\sum_{j\neq 1}v_j\sum_{l\neq0}{\partial U_{1,l}\over\partial x_1}dx=o(e^{-2\sigma_1}\sigma_1^{1-N\over2}).$$

Finally
\begin{equation}
\int_{\Omega_1}U_1^{p-1}(\sum_{j\neq1}v_j+\sum_{l\neq0}U_{1,l}){\partial v_1\over\partial x_1}dx=\int_{\Omega_1}U_1^{p-1}\sum_{j\neq1}v_j{\partial U_1\over\partial x_1}dx+o(e^{-2\underline\sigma}\underline\sigma^{1-N\over2}).
\end{equation}

Now  (\ref{eq:di2}), (\ref{eq:1er}), (\ref{eq:reste1}), (\ref{eq:jnon1}) and (\ref{eq:un}) give  the proof of the proposition.\\

{\bf Proof of Proposition \ref{unicite}. The points are asymptotically uniformly distributed.}\\

Let $i_0\in\{1,...,k\}$ be such that 
$$\sigma_{i_0}-\underline\sigma\rightarrow 0$$
(we know that there exists at least one $i$ such that $\sigma_{i}=\underline\sigma$).

By (\ref{eq:deq0}), we know that
$$d_{i_0}=0.$$

We deduce from (\ref{eq:di})   that
$$\int_{\Omega_{i_0}}U_{i_0}^{p-1}\sum_{j\neq i_0}v_j{\partial U_{i_0}\over\partial x_1}dx=o(e^{-2\underline\sigma}\underline\sigma^{1-N\over2}),$$

that is

$$\int_{\Omega_{i_0};x_1>\frac{a_{\ep}^{i_0}}{\ep}}U_{i_0}^{p-1}\sum_{j\neq i_0}v_j{\partial U_{i_0}\over\partial x_1}dx$$
$$=-\int_{\Omega_{i_0};x_1<\frac{a_{\ep}^{i_0}}{\ep}}U_{i_0}^{p-1}\sum_{j\neq i_0}v_j{\partial U_{i_0}\over\partial x_1}dx+o(e^{-2\underline\sigma}\underline\sigma^{1-N\over2}).
$$
Since we have
$$2\underline\sigma-(\frac{a_{\ep}^{i_0+1}}{\ep}-\frac{a_{\ep}^{i_0}}{\ep})\rightarrow0,$$
we use Corollary \ref{equivalence} to get a positive real number $D_0$ such that
$$e^{2\underline\sigma}\underline\sigma^{N-1\over2}\int_{\Omega_{i_0};x_1>\frac{a_{\ep}^{i_0}}{\ep}}U_{i_0}^{p-1}\sum_{j\neq i_0}v_j{\partial U_{i_0}\over\partial x_1}dx\rightarrow D_0$$
and
$$-e^{\vert\frac{a_{\ep}^{i_0}-a_{\ep}^{i_0-1}}{\ep}\vert}\vert\frac{a_{\ep}^{i_0}-a_{\ep}^{i_0-1}}{\ep}\vert^{N-1\over2}\int_{\Omega_{i_0};x_1<\frac{a_{\ep}^{i_0}}{\ep}}U_{i_0}^{p-1}\sum_{j\neq i_0}v_j{\partial U_{i_0}\over\partial x_1}dx\rightarrow D_0.$$
and consequently

$$2\underline\sigma-(\frac{a_{\ep}^{i_0}}{\ep}-\frac{a_{\ep}^{i_0-1}}{\ep})\rightarrow 0.$$
This property is valid for all exponent $i$ such that $\sigma_i-\underline\sigma\rightarrow 0$ instead of $i_0$. Thus
we have (\ref{eq:equirepartis}).\\

\section{ The proof of Theorem \ref{uniqueDancer} completed.}

{\bf { The uniqueness.}}

Now, we have $\frac{a_{\ep}^i}{\ep}-(\frac{a_{\ep}^1}{\ep}+\frac{(i-1)2\pi}{k\ep})\rightarrow0.$ Replacing the points $\frac{a_{\ep}^i}{\ep}$ by $\frac{a_{\ep}^1}{\ep}+\frac{(i-1)2\pi}{k\ep}$, $i=1,...,k$, we can write $u$ as
$$u=\sum_{l\in\Z}U(x_1-\frac{a_{\ep}^1}{\ep}+\frac{2\pi l}{k\ep},x')+\sum_{i=1}^k\tilde\delta_i\tilde\varphi_i+\tilde v,\quad<\tilde v,\tilde\varphi_i>=0,\quad i=1,...,k.$$
By the definition of $\tilde\varphi_i$ given in section 2 (analogue to that of $\varphi_i$), $\tilde\varphi_i$ is $\frac{2\pi}{\ep}$-periodic in $x_1$. But now, 
the corresponding operator $\LL$ is of minimal period $2\pi\over k\ep$, since now $\overline u_{\ep}$ is replaced by $\sum_{l\in\Z}U(x_1-\frac{a_{\ep}^1}{\ep}+\frac{2\pi l}{k\ep},x')$. So we have
$$\tilde\varphi_i(x)=\tilde\varphi_1(x_1+\frac{2(i-1)\pi}{k\ep},x')$$
and $\tilde\varphi_1$ is $\frac{2\pi}{k\ep}$-periodic. Let us denote $\varphi_1=\tilde\varphi$.\\
Now we recall that
$$\LL \tilde v=h^{\bot}$$
with $h=-\mathcal{M}(\overline u_{\ep})+O(\tilde v^2+(\sum_i\tilde\delta_i\tilde\varphi_i)^2)$.
We can use the Banach fixed point Theorem in $L^{\infty}$ to deduce that $\tilde v$ is of minimal period $2\pi\over k\ep$.\\
 Consequently, $u$ is $2\pi\over k\ep$-periodic and in the space $H^1(\frac{S^1}{k\ep}\times\R^{N-1})$ we write
$$u=\sum_{l\in\Z}U(x_1-\frac{a_{\ep}^1}{\ep}+\frac{2\pi l}{k\ep},x')+\tilde\delta\tilde\varphi+\tilde v,\quad<\tilde v,\tilde\varphi>=0.$$
Thanks to Proposition \ref{deltai}, we can perform a translation in $x_1$ to get $\tilde\delta=0.$ We get some $\frac{a_{\ep}}{\ep}\rightarrow 0$ such that
$$u=\sum_{l\in\Z}U(x_1-\frac{a_{\ep}^1}{\ep}-\frac{a_{\ep}}{\ep}+\frac{2\pi l}{k\ep},x')+v,$$
\begin{equation}\label{eq:orthphi}
< v,\tilde\varphi(x_1-\frac{a_{\ep}}{\ep})>=0.
\end{equation}
Let $u_1$ be another solution of period $2\pi\over \ep$, which verify (\ref{eq:conv}).
Exactly as for $u$, we find some point $\frac{b_{\ep}}{\ep}\rightarrow0$ such that
$$u_1=\sum_{l\in\Z}U(x_1-\frac{a_{\ep}^1}{\ep}-\frac{b_{\ep}}{\ep}+\frac{2\pi l}{k\ep},x')+v_1,\quad<v_1,\tilde\varphi(x_1-\frac{b_{\ep}}{\ep})>=0.$$
We define
$$\overline u_1(x)=u_1(x_1+\frac{b_{\ep}-a_{\ep}}{\ep},x')\quad \hbox{and}\quad 
w=u- \overline u_1.$$
Now we can prove that for $\ep$ small enough, 
$$w=0.$$
The proof is very close to the proof of the uniqueness of the Dancer solution among the even solutions verifying
(\ref{eq:Dancer}), but here the property of being even in $x_1$ is replaced by the property (\ref{eq:orthphi}). Let us give the proof for the sake of completeness.\\
Without loss of generality, we can suppose that $a_{\ep}^1=0$. 
Let us suppose that $w\neq0$, at least for a sequence $\ep\rightarrow0$. Then $\|w\|_{\infty}$ is attained at a point $c=(c_1,c')$, with $c'$ obviously bounded independently of $\ep$ and $c_1\in]-\frac{\pi}{k\ep},\frac{\pi}{k\ep}]$. Now $c_1$ is bounded. To see that, we write
\begin{equation}\label{eq:unicite}
-\Delta w+ w(1-\frac{u^p-v^p}{u-v})=0.
\end{equation}
If $w(c)>0$, then 
$$\frac{u^p-\overline u^p_1}{u-\overline u_1}(c)\leq pu^{p-1}$$
thus
$$\Delta w(c) \geq w(c)(1- pu^{p-1}(c)).$$
But if $\vert c_1\vert\rightarrow+\infty$, we have $u^{p-1}(c)\rightarrow0$, thus
$$1- pu^{p-1}(c)>0$$
for $\ep$ small enough
that is in contradiction with the Maximum Principle.
So we may extract a subsequence such that $c\rightarrow\overline c$ for some $\overline c$.\\
Let us define
$$z(x) =\frac{w}{\|w\|_{\infty}}.$$
It verifies (\ref{eq:unicite}). By standard arguments $z\rightarrow\overline z$ uniformly on compact sets. Moreover $$\overline z(\overline c)=1$$ and
$$p\overline u_1^{p-1}\leq\frac{u^p-\overline u_1^p}{u-\overline u_1}\leq pu^{p-1}\quad\hbox{if $u>\overline u_1$}$$
and we have the reverse inequality if $u<\overline u_1$. More,
$$\lim u=\lim \overline u_1=U(x),$$
uniformly on compact sets.
So
$$\lim\frac{u^p-\overline u_1^p}{u-\overline u_1}=pU^{p-1}.$$
Thus
$$-\Delta\overline z+\overline z(1-pU^{p-1})=0.$$
We deduce that $\overline z=\alpha{\partial U\over\partial x_1}$, for some $\alpha\neq0$.

We have
$$<z,\tilde\varphi(x_1-\frac{a_{\ep}}{\ep},x')>_{H^1}={1\over\|u-\overline u_1\|_{\infty}}<v-v_1(x_1+\frac{b_{\ep}-a_{\ep}}{\ep},x'),\tilde\varphi(x_1-\frac{a_{\ep}}{\ep},x')>_{H^1}=0.$$
 Moreover
$$\vert z\vert\leq 1,\quad \tilde\varphi\rightarrow{\partial U\over\partial x_1},$$
and by Proposition \ref{barrierefonctionpropre},
$$\vert\nabla\tilde\varphi(x)\vert+\vert\tilde\varphi(x)\vert\leq Ce^{-{\eta\vert x\vert}}\quad\hbox{in}\quad[-\frac{\pi}{k\ep},\frac{\pi}{k\ep}].$$
We use the Lebesgue Theorem to infer that
$$<z,\tilde\varphi(x_1-\frac{a_{\ep}}{\ep},x')>_{H^1}\rightarrow \alpha\|{\partial U\over\partial x_1}\|^2_{H^1(\R^N)}.$$

 So we are led to a contradiction. We conclude that $u=\overline u_1$, for $\ep$ small enough.
 In particular, we may take the solution $u_D(x_1-\frac{a_{\ep}^1}{\ep},x')$ instead of $u_1$.
 \\

{\bf Estimating the difference between the Dancer solution and the groundstate solution. (Proof of \ref{eq:estimationpsi})}.\\

Without loss of generality, we let $k=1$. We write

$$u_D=\sum_lU_l+v$$
where $v$ is even in $x_1$ and verifies 
$$\LL v=h$$
and $\overline u_{\ep}$ is replaced by $\sum_lU_l$ in the definition of $\LL$. The restriction of $\LL$ to the even functions has no eigenvalue tending to 0. A proof very close to the proof of (\ref{eq:xiinfini}) gives
$$\|v\|_{\infty}+\|\nabla v\|_{\infty}\leq C\|h\|_{\infty}$$
and consequently, as for (\ref{eq:normevinfini}) (with $\delta_i=0$), we deduce
\begin{equation}\label{eq:estiv}
 \quad\|v\|_{\infty}+\|\nabla v\|_{\infty}\leq C e^{-{2\pi\over\ep}}({\pi\over\ep})^{1-N\over2}.      
\end{equation}

Let $R_0>0$ be given. It remains to estimate, for all $\eta\in]0,1[$
$$(\vert v(y)\vert+\vert\nabla v(y)\vert)e^{\eta d_y}$$
when $d_y>R_0,$
where we denote
$$d_y=\hbox{dist}(y,\cup_{l\in\Z}\{(\frac{2l\pi}{\ep},0)\}).$$

 Let
 $$y\in ]-{\pi\over\ep},{\pi\over\ep}[\times\R^{N-1}.$$ 
 We follow the course of the proof of Proposition \ref{appendix1} from (\ref{eq:theta}), $\xi$ being replaced by $v$.
With the notations of that proof, we let $\tilde\sigma_i=\tilde\sigma_{i-1}=\frac{\pi}{\ep}$ and we consider a positive real number $\beta$, independent of $\ep$, which will be chosen later. We perform the truncation around 0, using the truncation function $\theta$. So we drop the index $i$. By (\ref{eq:Greenbis}), we have
$$(\theta v)(y)=\int_{\R^{N}} G(y-x)(p\theta\overline u_{\ep}^{p-1} v+\theta h-\Delta\theta v-2\nabla\theta.\nabla v)(x)dx$$
Now $$h=-\mathcal{M}(\overline u_{\ep})+O(v^2).$$

Since 
$\theta(y)=1$, we have
$$\vert v(y)\vert\leq C(
\int_{\R^{N}} vG(y-x)p\theta e^{(-p+1)d_x}dx+\int_{\R^{N}} G(y-x)\theta\vert\mathcal{M}(\overline u_{\ep})\vert dx
+\int_{\R^{N}} G(y-x)\theta v^2dx$$
$$
+\frac{1}{\beta}\int_{\R^{N};\frac{\pi}{\ep}<x_1<\frac{\pi}{\ep}+\beta} G(y-x)(\vert v\vert+\vert\nabla v\vert)dx).$$

 First we write, for all $0<\eta<1$
$$\int_{\R^{N}} \vert v\vert G(y-x)p\theta e^{(-p+1)d_x}dx\leq \|ve^{\eta d_x}\|_{\infty}\int_{\R^{N}}  G(y-x)p\theta e^{(-p+1-\eta)d_x}dx.
$$
For $x\in$Supp$\theta$ we have
$$d_x=\vert x-(\frac{l2\pi}{\ep},0)\vert,\quad \hbox{$l=0$, or $l=1$, or $l=-1$}.$$
Then
$$\vert x-y\vert\geq \vert y-(\frac{l2\pi}{\ep},0)\vert-\vert x-(\frac{l2\pi}{\ep},0)\vert\geq d_y-d_x.$$
In view of (\ref{eq:G}), we get
$$G(y-x)\leq e^{-d_y+d_x}.$$
Then
\begin{equation}\label{eq:1}
\int_{\R^{N}}  G(y-x)p\theta e^{(-p+1-\eta)d_x}dx\leq C e^{-d_y }\int_{\R^{N}}  p\theta e^{(-p+2-\eta)d_x}\leq  Ce^{-\eta d_y} e^{-\eta_1 R_0 },
\end{equation}
where we define $\eta_1=1-\eta$.\\
Now, we have
 $$\quad U\leq Ce^{-\vert x\vert} \vert x\vert^{1-N\over2}$$
and in $]-\frac{\pi}{\ep},\frac{\pi}{\ep}[\times \R^{N-1}$
\begin{equation}\label{eq:M3}
\mathcal{M}(\overline u_{\ep})=pU^{p-1}(\sum_{l\neq0}U_l)+O((\sum_{l\neq0}U_l)^2.
\end{equation}

We deduce that for all $\eta\in]0,1[$

$$
\int_{]-\frac{\pi}{\ep},\frac{\pi}{\ep}[\times \R^{N-1}} G(y-x)\theta\vert\mathcal{M}(\overline u_{\ep})\vert dx\leq C\sum_{l\in\Z^{\star}}\int_{]-\frac{\pi}{\ep},\frac{\pi}{\ep}[\times \R^{N-1}}e^{-\eta\vert y\vert} e^{\eta\vert x-y\vert}G(y-x)e^{\eta\vert x\vert}(U^{p-1} U_l+U_l^2)\theta dx.
$$
  We obtain
$$\sum_{l\in\Z^{\star}}\int_{]-\frac{\pi}{\ep},\frac{\pi}{\ep}[\times \R^{N-1}}e^{-\eta\vert y\vert} e^{\eta\vert x-y\vert}G(y-x)e^{\eta\vert x\vert}U^{p-1} U_l\theta dx$$
$$\leq C\sum_{l\in\Z^{\star}}\int_{]-\frac{\pi}{\ep},\frac{\pi}{\ep}[\times \R^{N-1}}e^{-\eta d_y} e^{(-1+\eta)\vert x-y\vert}e^{(-p+1+\eta)\vert x\vert}e^{-\vert x-\frac{2l\pi}{\ep}\vert}\vert x-\frac{2l\pi}{\ep}\vert^{1-N\over2}\theta dx$$
We have
$$\vert x-(\frac{2l\pi}{\ep},0)\vert\geq \vert \frac{2l\pi}{\ep}\vert-\vert x\vert\geq \frac{2\pi}{\ep}-\vert x\vert$$
and, when $x\in]-\frac{\pi}{\ep},\frac{\pi}{\ep}[\times \R^{N-1}$ and $l\in\Z^{\star}$
$$\vert x-(\frac{2\pi l}{\ep},0)\vert
\geq\frac{\pi}{\ep}.$$
We deduce that for all $0<\eta'<1$ such that $\eta'<p-1-\eta$
$$\sum_{l\in\Z^{\star}}\int_{]-\frac{\pi}{\ep},\frac{\pi}{\ep}[\times \R^{N-1}}e^{-\eta d_y} e^{(-1+\eta)\vert x-y\vert}e^{(-p+1+\eta)\vert x\vert}e^{-\vert x-\frac{2l\pi}{\ep}\vert}\vert x-\frac{2l\pi}{\ep}\vert^{1-N\over2}\theta dx$$
$$\leq Ce^{-\eta d_y}e^{-\eta'\frac{2\pi}{\ep}}(\frac{\pi}{\ep})^{1-N\over2}\int_{]-\frac{\pi}{\ep},\frac{\pi}{\ep}[\times \R^{N-1}}\theta e^{(-p+1+\eta+\eta')\vert x\vert}dx\leq Ce^{-\eta d_y}e^{-\eta'\frac{2\pi}{\ep}}(\frac{\pi}{\ep})^{1-N\over2}.$$
Now, when $x\in$Supp$\theta\cap]\frac{\pi}{\ep},\frac{\pi}{\ep}+\beta[\times \R^{N-1} $ and $l\in\Z^{\star}$, we have
\begin{equation}\label{eq:M4}
\mathcal{M}(\overline u_{\ep})=pU_1^{p-1}(\sum_{l\neq1}U_l)+O((\sum_{l\neq1}U_l)^2).
\end{equation}
Then
$$
\int_{]\frac{\pi}{\ep},\frac{\pi}{\ep}+\beta[\times \R^{N-1}} G(y-x)\theta\vert\mathcal{M}(\overline u_{\ep})\vert dx$$
$$\leq C\sum_{l\neq1}\int_{]\frac{\pi}{\ep},\frac{\pi}{\ep}+\beta[\times \R^{N-1}}e^{-\eta\vert y-\frac{2\pi}{\ep}\vert} e^{\eta\vert x-y\vert}G(y-x)e^{\eta\vert x-\frac{2\pi}{\ep}\vert}(U_1^{p-1} U_l+U_l^2)\theta dx.
$$
We have
$$\sum_{l\neq1}\int_{]\frac{\pi}{\ep},\frac{\pi}{\ep}+\beta[\times \R^{N-1}}e^{-\eta\vert y-\frac{2\pi}{\ep}\vert} e^{\eta\vert x-y\vert}G(y-x)e^{\eta\vert x-\frac{2\pi}{\ep}\vert}U_1^{p-1} U_l\theta dx$$
$$\leq C\int_{]\frac{\pi}{\ep},\frac{\pi}{\ep}+\beta[\times \R^{N-1}}e^{-\eta d_y} e^{(-1+\eta)\vert x-y\vert}e^{(-p+1+\eta)\vert x-\frac{2\pi}{\ep}\vert}e^{-\vert x\vert}\vert x\vert^{1-N\over2}\theta dx.$$
For $x\in]\frac{\pi}{\ep},\frac{\pi}{\ep}+\beta[\times \R^{N-1}$, we have
$$\vert x\vert\geq x_1\geq\frac{\pi}{\ep}$$
then, for all $0<\eta'<1$ such that $\eta'<p-1-\eta$
$$\int_{]\frac{\pi}{\ep},\frac{\pi}{\ep}+\beta[\times \R^{N-1}}e^{-\eta d_y} e^{(-1+\eta)\vert x-y\vert}e^{(-p+1+\eta)\vert x-\frac{2\pi}{\ep}\vert}e^{-\vert x\vert}\vert x\vert^{1-N\over2}\theta dx$$
$$\leq C\int_{]\frac{\pi}{\ep},\frac{\pi}{\ep}+\beta[\times \R^{N-1}}e^{-\eta d_y} e^{(-p+1+\eta+\eta')\vert x-\frac{2\pi}{\ep}\vert}e^{-\eta'2\frac{\pi}{\ep}}(\frac{\pi}{\ep})^{1-N\over2}\leq Ce^{-\eta d_y}e^{-\eta'\frac{2\pi}{\ep}}(\frac{\pi}{\ep})^{1-N\over2}.$$
Now, we have
 $$\vert x-(\frac{2\pi }{\ep},0)\vert\geq\frac{\pi}{\ep}\geq\vert x\vert,\quad\hbox{for $x\in]-\frac{\pi}{\ep},\frac{\pi}{\ep}[\times\R^{N-1}$}.$$
We deduce that
$$\sum_{l\in\Z^{\star}}\int_{]-\frac{\pi}{\ep},\frac{\pi}{\ep}[\times\R^{N-1}}e^{-\eta\vert y\vert} e^{\eta\vert x-y\vert}G(y-x)e^{\eta\vert x\vert} U_l^2\theta dx$$
$$\leq C\int_{]-\frac{\pi}{\ep},\frac{\pi}{\ep}[\times\R^{N-1}}e^{-\eta d_y}e^{(-1+\eta)\vert x-y\vert}e^{(-2+\eta)\vert x-(\frac{2\pi }{\ep},0)\vert}dx(\frac{\pi}{\ep})^{1-N}$$
$$\leq Ce^{-\eta d_y}\int_{]-\frac{\pi}{\ep},\frac{\pi}{\ep}[\times\R^{N-1}}e^{(-1+\eta)\vert x-y\vert}dxe^{(-2+\eta)\frac{\pi }{\ep}}(\frac{\pi}{\ep})^{1-N}\leq Ce^{-\eta d_y}e^{(-2+\eta)\frac{\pi }{\ep}}(\frac{\pi}{\ep})^{1-N}.
$$

In the same way, we get

$$\sum_{l\neq1}\int_{]\frac{\pi}{\ep}+\beta,\frac{\pi}{\ep}[\times\R^{N-1}}e^{-\eta\vert y-\frac{2\pi}{\ep}\vert} e^{\eta\vert x-y\vert}G(y-x)e^{\eta\vert x-\frac{2\pi}{\ep}\vert} U^2\theta dx$$
$$\leq C\int_{]\frac{\pi}{\ep}+\beta,\frac{\pi}{\ep}[\times\R^{N-1}}e^{-\eta d_y}e^{(-1+\eta)\vert x-y\vert}e^{(-2+\eta)\vert x\vert}dx(\frac{\pi}{\ep})^{1-N}$$
$$\leq Ce^{-\eta d_y}\int_{]-\frac{\pi}{\ep},\frac{\pi}{\ep}[\times\R^{N-1}}e^{(-1+\eta)\vert x-y\vert}e^{(-2+\eta)\frac{\pi }{\ep}}dx(\frac{\pi}{\ep})^{1-N}\leq Ce^{-\eta d_y}e^{(-2+\eta)\frac{\pi }{\ep}}(\frac{\pi}{\ep})^{1-N}.
$$
We have proved that, for all $0<\eta'<p-1-\eta$ and $\eta'<1$

$$
\int_{\R^{N}} G(y-x)\theta\vert\mathcal{M}(\overline u_{\ep})\vert dx\leq Ce^{-\eta d_y}(\frac{\pi}{\ep})^{1-N\over2}(e^{-\eta'\frac{2\pi}{\ep}}+e^{(-2+\eta)\frac{\pi}{\ep}}(\frac{\pi}{\ep})^{1-N\over2}).
$$

Then
, we obtain, for all $0<\eta'<p-1-\eta$ and $0<\eta'<1$
\begin{equation}\label{eq:2}
\int_{\R^{N}} G(y-x)\theta\vert\mathcal{M}(\overline u_{\ep})\vert dx\leq Ce^{-\eta d_y}(\frac{\pi}{\ep})^{1-N\over2}(e^{-2\eta'\frac{\pi}{\ep}})
\end{equation}

Now, since $d_x\geq d_y-\vert x-y\vert$, we have

\begin{equation}\label{eq:3}
\int_{\R^{N}} G(y-x)v^2\theta dx\leq\|ve^{\eta d_x}\|^2_{\infty}\int_{\R^{N}} G(y-x)\theta e^{-2\eta d_x}dx
\end{equation}
$$\leq C\|ve^{\eta d_x}\|^2_{\infty}e^{-\eta d_y}\int_{\R^{N}} G(y-x)\theta e^{\eta\vert x-y\vert}e^{-\eta d_x}dx\leq C\|ve^{\eta d_x}\|^2_{\infty}e^{-\eta d_y}.$$

Last
\begin{equation}\label{eq:4}
\int_{\R^{N};\frac{\pi}{\ep}<x_1<\frac{\pi}{\ep}+\beta} G(y-x)(\vert v\vert+\vert\nabla v\vert)dx
\leq(\|ve^{\eta d_x}\|_{\infty}
+\|(\nabla v)e^{\eta d_x} \|_{\infty})e^{-\eta d_y}\int_{\R^{N}}e^{\eta\vert x-y\vert} G(y-x) dx.
\end{equation}
Finally, (\ref{eq:1}), (\ref{eq:2}), (\ref{eq:3}) and (\ref{eq:4}) give for all $y\in ]-{\pi\over\ep},{\pi\over\ep}[\times\R^{N-1}$ such that $d_y\geq R_0$ and for all  $0<\eta'<1$, $0<\eta'<p-1-\eta$ and for $\eta_1=1-\eta$
$$|v(y)e^{\eta \vert y\vert}|\leq Ce^{-{2\eta'\pi\over\ep}}({\pi\over\ep})^{1-N\over2}+{C\over\beta}(\|ve^{\eta d_x}\|_{\infty}+\|(\nabla v)e^{\eta d_x}\|_{\infty})+C(\|ve^{\eta d_x}\|^2_{\infty}+\|ve^{\eta d_x}\|_{\infty} e^{-\eta_1 R_0 }).
$$
We already know that if $d_y\leq R_0$, then
$$(\vert\nabla v(y)\vert+\vert v(y)\vert)e^{\eta d_y} \leq C e^{-{2\eta'\pi\over\ep}}({\pi\over\ep})^{1-N\over2}e^{\eta R_0}.$$
Thus 
$$|v(y)e^{\eta \vert y\vert}|\leq Ce^{-{2\eta'\pi\over\ep}}e^{\eta R_0}({\pi\over\ep})^{1-N\over2}+{C\over\beta}(\|ve^{\eta d_x}\|_{\infty}+\|(\nabla v)e^{\eta d_x}\|_{\infty})+C(\|ve^{\eta d_x}\|^2_{\infty}+\|ve^{\eta d_x}\|_{\infty} e^{-\eta_1 R_0 })
$$
for all $y$, where the constants are independent of $\beta$ and of $R_0$.\\
We obtain the same estimate for $|(\nabla v)(y)e^{\eta \vert y\vert}|$, using the proof of (\ref{eq:h1suite}). \\
Now, we use (\ref{eq:estiv}) to obtain for all $x$
$$(\vert v(x)\vert+\vert\nabla v(x)\vert)e^{d_x}\leq(\vert v(x)\vert+\vert\nabla v(x)\vert)e^{\frac{\pi}{\ep}}\leq
e^{-{\pi\over\ep}}\rightarrow0.$$

Thus, we can choose $\beta$ and $R_0$ large enough and $\ep$ small enough to obtain (\ref{eq:estimationpsi}).\\

\section{Appendix.}
 
Let $\eta>0$ be given. Let $h$ be a function, defined on $[-\frac{\pi}{\ep},\frac{\pi}{\ep}]\times\R^{N-1}$, which has the following property 
\begin{equation}\label{eq:proprh1}
\hbox{ $x\mapsto h(x)e^{\eta\hbox{dist}(x,\cup_{j=0}^{k+1}\{(\frac{a_{\ep}^j}{\ep},0)\})}$ 
is bounded in $ L^{\infty}([-\frac{\pi}{\ep},\frac{\pi}{\ep}]\times\R^{N-1})$}
\end{equation}
independently of $\ep$.\\

Then $h$ belongs to $H^{-1}([-\frac{\pi}{\ep},\frac{\pi}{\ep}]\times\R^{N-1})$(the dual space of $H^{1}([-\frac{\pi}{\ep},\frac{\pi}{\ep}]\times\R^{N-1})$) in the following sense
$$<h,\psi>_{H^{-1},H^1}=\int_{[-\frac{\pi}{\ep},\frac{\pi}{\ep}]\times\R^{N-1}}h\psi dx\quad\hbox{for $\psi\in H^1$}.$$
We will denote $H^1$ in place of $H^{1}([-\frac{\pi}{\ep},\frac{\pi}{\ep}]\times\R^{N-1})$.\\
By the Lax-Milgram Theorem, there exists $u\in H^1$  such that
\begin{equation}\label{eq:h-1}
-\Delta u+u=h.
\end{equation}
It is classical that $u\in L^{\infty}$ (see \cite{GT}, Theorem 9.13 and use the Sobolev embedding Theorem). 

As a consequence of the maximum principle, 
\begin{equation}\label{eq:PM}
\|u\|_{L^{\infty}}\leq \| h\|_{L^{\infty}}. 
\end{equation}
More we have the following 

\begin{proposition}\label{appendix1}

Let $\eta\in]0,1[$ be given. Let us suppose that $h\in L^{\infty}([-\frac{\pi}{\ep},\frac{\pi}{\ep}]\times\R^{N-1})$ has the property (\ref{eq:proprh1}) and that
 \begin{equation}\label{eq:proprh2}
<h,\varphi_i>_{H^{-1},H^1}=0,\quad i=1,...,k.
\end{equation}
  Then  there exists a unique $\xi \in H^1([-\frac{\pi}{\ep},\frac{\pi}{\ep}]\times\R^{N-1})$ which verifies
\begin{equation}\label{eq:Lh}
\LL\xi= h.
\end{equation}
\begin{equation}\label{eq:xiinfini}
\|\xi \|_{\infty}+\|\nabla\xi\|_{\infty}\leq C \| h \|_{\infty}
\end{equation}
where $C$ is independent of $\ep$.\\

More, there exists $C$ independent of $\ep$ and dependent of $\eta$ such that
 \begin{equation}\label{eq:hj1}
 \| \xi e^{\eta \hbox{dist}(x,\cup_{j=0}^{k+1}\{(\frac{a_{\ep}^j}{\ep},0)\})}\|_{\infty}+\| \nabla\xi e^{\eta \hbox{dist}(x,\cup_{j=0}^{k+1}\{(\frac{a_{\ep}^j}{\ep},0)\})}\|_{\infty}
 \end{equation}
 $$\leq C\|he^{\eta\hbox{dist}(x,\cup_{j=0}^{k+1}\{(\frac{a_{\ep}^j}{\ep},0)\})}\|_{\infty}.$$
\end{proposition}

{\bf Proof} 
The operator $\LL$ being a Fredholm operator, we have the existence of a unique solution $\xi$ of (\ref{eq:Lh}) when $h\in H^{-1}$ and verifies the property (\ref{eq:proprh2}). In this case, $\xi$ verifies
\begin{equation}\label{eq:orth}
\int_{{S^1\over\ep}\times\R^{N-1}}\xi(-\Delta+1)(\varphi_i) dx=0,\quad\hbox{i=1,...,k}.
\end{equation}
Moreover we have in this case
\begin{equation}\label{eq:minorant}
\|\xi\|_{H^1}\leq C\| u\|_{H^1}
\end{equation}
where $u$ is defined in (\ref{eq:h-1}) and $C$ is independent of $\ep$. Indeed, this can be proved using the expansion of $\xi$ on a basis of eigenvectors of the operator $(-\Delta+1)^{-1}\LL$.\\
 Now, let $h\in L^{\infty}$ verifying (\ref{eq:proprh1}) and (\ref{eq:proprh2}). Let us prove (\ref{eq:xiinfini}).\\
By a standard proof, we have that $\xi$ is continuous. \\
Indeed, $\xi\in H^1$. By the Sobolev embedding theorem (see \cite{GT}, Theorem 7.26) we get $\xi\in L^q$, with $q=\frac{2N}{N-2}$ if $N>2$ and $\xi\in L^q$ for all $q\geq2$ if $N=2$. By (\ref{eq:proprh1}), $h\in L^r$ for all $1\leq r\leq\infty$. We get that $-\Delta\xi\in L^q$. This gives that $\xi\in W^{2,q}$ (see \cite{GT}, Theorem 9.9). We use the Sobolev embedding theorem again. If $N=2$, we choose $q>{N\over2}$ to get that $\xi$ is continuous. If $N>2$, since $\frac{N}{N-2}>1$, iterating the process, we will reach $q>{N\over2}$ such that $\xi\in W^{2,q}$ and we will conclude that $\xi$ is continuous.\\
 Let us assume that $\|h\|_{\infty}\rightarrow0$ and that $\|\xi\|_{\infty}=1$. Let us remark that, using the periodic Green function as in \cite{estimate}, we can prove very quickly the following provisional estimate 
$$\exists C\quad
\vert\xi(x_1,x')\vert\leq Ce^{-\eta\vert x'\vert}.
$$
Now, since $\xi$ is continuous and since $\xi(x_1,x')\rightarrow0$ as $\vert x'\vert\rightarrow0$, there exists $c$ such that $\xi(c)=1$. Let $\tilde\xi(x)=\xi(x+c)$. By standard elliptic estimates, $\tilde\xi$ tends to a limit $\overline\xi$, uniformly on compact sets of $\R^N$ and we have either 
$$(-\Delta+1)\overline\xi=0\quad\hbox{ in $\R^N$ if $\vert c-\frac{a_{\ep}^i}{\ep}\vert\rightarrow+\infty$ for all $i$}$$
or
$$(-\Delta+1-pU^{p-1}(x+\overline c))\overline\xi=0\quad\hbox{ in $\R^N$ if there exists $i$ and $\overline c$ such that $(c-\frac{a_{\ep}^i}{\ep})\rightarrow\overline c$. }$$
The first case is in contradiction with the maximum principle, so it does not occur. In the second case, we have that
$$\xi(x+c)\rightarrow\frac{\partial U}{\partial x_1}(x_1+\overline c,x')\quad\hbox{uniformly on compact sets}.$$
We use (\ref{eq:orth}). Since $(1-\lambda_i)(-\Delta+1)\varphi_i=p\overline u_{\ep}^{p-1}\varphi_i$, we use the Lebesgue Theorem to get a contradiction.\\

We have proved that
$$\|\xi\|_{\infty}\leq C\|h\|_{\infty}.$$
The inequality for $\|\nabla\xi\|_{\infty}$ follows from standard elliptic estimates (\cite{GT},Theorem 9.13). So we have proved (\ref{eq:xiinfini}).\\

In what follows, we denote
$$d_x=\hbox{dist}(x,\cup_{j=0}^{k+1}\{(\frac{a_{\ep}^j}{\ep},0)\}).$$

 We define
$$\tilde \sigma_i=\frac{a_{\ep}^{i+1}-a_{\ep}^{i}}{2\ep}\quad\hbox{for $i=0,...,k$}.$$
Let $y\in]-\frac{\pi}{\ep},\frac{\pi}{\ep}[\times\R^{N-1}$ be given. Then
$$(y\in\Omega_i)\Leftrightarrow
(\tilde\sigma_{i-1}\leq y_1-\frac{a_{\ep}^i}{\ep}\leq \tilde\sigma_{i}),$$
where $\Omega_i$ is defined in (\ref{eq:defomega}).\\
Let $\beta$ be a positive real number, independent of $\ep$. 
The number $\beta$ will be chosen later. \\
Let $R_0>0$ be a given real number, independent of $\ep$. We are going to estimate
$$\vert \xi(y)e^{\eta d_y}\vert+\vert \nabla\xi(y)e^{\eta d_y}\vert$$
when $y\in\Omega_i$ is such that $d_y>R_0$.\\

For $i=1,...,k$, let $\theta_i$ be a function which verifies

\begin{equation}\label{eq:theta}
\theta_i(x)=1\quad\hbox{for}\quad x\in\Omega_i
\end{equation}
$$\theta_i(x)=0\quad\hbox{for}\quad \quad x_1-\frac{a_{\ep}^i}{\ep}\geq\tilde\sigma_{i}+\beta\quad\hbox{or}\quad x_1-\frac{a_{\ep}^i}{\ep}\leq -\tilde\sigma_{i-1}-\beta.$$
Moreover, we suppose that $\theta_i$ is $\mathcal{C}^2$ in $x_1$. More precisely, we build $\theta_i$ from the function $\tilde\theta$ defined in $[0,\beta]$, $\tilde\theta(x_1)=-{6\over\beta^5}x_1^5+{15\over\beta^4}x_1^4-{10\over\beta^3}x_1^3+1$, by $\theta_i(x)=\tilde\theta(x_1-{a_{\ep}^i\over\ep}-\tilde\sigma_i)$, if $\tilde\sigma_{i}\leq x_1-\frac{a_{\ep}^i}{\ep}\leq\tilde\sigma_{i}+\beta$ and $\theta_i(x)=\tilde\theta( {a_{\ep}^i\over\ep}-\tilde\sigma_{i-1}-x_1)$, if $-\tilde\sigma_{i-1}-\beta\leq x_1-\frac{a_{\ep}^i}{\ep}\leq-\tilde\sigma_{i-1}.$
\\
Thus we have for all $x$ and for $M$ independent of $i$, of $\beta$ and of $\ep$
$$\vert\theta_i(x)\vert\leq M,\quad\vert\nabla\theta_i(x)\vert\leq {M\over\beta},\quad\vert\Delta\theta_i(x)\vert\leq {M\over\beta^2}.$$
 Let $G$ be the Green function of the operator 
$$-\Delta+1\quad\hbox{ on }\R^N$$
with the null Dirichlet condition at infinity. It is recalled in \cite{Gidas} that
\begin{equation}\label{eq:G}
0<G(r)\leq C{e^{-r}\over r^{N-2}}(1+r)^{(N-3)/2}\hbox{  for  } N\geq 2
\end{equation}
and we have the following estimate, valid for all $\eta\in]0,1[$.
\begin{equation}\label{eq:integraleeta}
\int_{\R^{N}}G( y-x)e^{-\eta\vert y\vert}dy\leq Ce^{-\eta\vert x\vert}
\end{equation}
which is an easy consequence of (\ref{eq:G}).
We write
 \begin{equation}\label{eq:Greenbis}
 (\theta_i \xi)(y)=\int_{\R^{N}} G(y-x)(p\theta_i\overline u_{\ep}^{p-1}\xi+\theta_ih-\Delta\theta_i\xi-2\nabla\theta_i.\nabla\xi)(x)dx,
\end{equation}

Since $y\in{\Omega_i}$, we have $(\theta_i\xi)(y)=\xi(y)$.\\
We consider (\ref{eq:Greenbis}).
Firstly, we have for all $\eta\in]0,1[$
\begin{equation}\label{eq:part0}
\vert\int_{\R^{N}} G(y-x)(p\theta_i\overline u_{\ep}^{p-1}\xi+\theta_ih)(x)dx\vert\leq C\|he^{\eta d_x}\|_{\infty}\int_{\R^N}G(y-x)\theta_i(x)e^{-\eta d_x}dx
\end{equation}
$$+C\|\xi e^{\eta d_x}\|_{\infty}\int_{\R^N}G(y-x)\theta_i(x)e^{-\eta d_x}\overline u_{\ep}^{p-1}(x)dx.
$$

We have for all $x$
$$d_x\geq d_y-\vert x-y\vert$$

So we have
\begin{equation}\label{eq:part2}
\vert\int_{\R^N} G(y-x)\theta_i(x) e^{-\eta d_x}dx\vert
\end{equation}
$$\leq\int_{\R^N} G(y-x)e^{\eta\vert y-x\vert}\vert\theta_i(x)\vert e^{-\eta d_y}dx\leq Ce^{-\eta d_y}.
$$
Now, for all $x$
$$\quad\overline u_{\ep}(x)\leq Ce^{- d_x}.$$
Then, for any $\eta_1$ such that $0<\eta_1<1-\eta$ 
\begin{equation}\label{eq:part3}\vert\int_{\R^N} G(y-x)\theta_i(x)e^{-\eta d_x}\overline u_{\ep}^{p-1}(x)dx\vert
\end{equation}
$$\leq Ce^{-(\eta+\eta_1)d_y}\int_{\R^N} G(y-x)e^{(\eta+\eta_1)\vert x-y\vert}e^{-(p-1-\eta_1)d_x}\vert\theta_i(x)\vert dx\leq Ce^{-(\eta+\eta_1)d_y}.$$

Now, we get for all $\eta\in]0,1[$
\begin{equation}\label{eq:part7}
\vert\int_{\R^{N}} G(y-x)(-\Delta\theta_i\xi-2\nabla\theta_i.\nabla\xi)(x)dx\vert\leq {C\over\beta}(\|\xi e^{\eta d_x}\|_{\infty}
\end{equation}
$$+\|\nabla\xi e^{\eta d_x}\|_{\infty})\int_{x\in\hbox{Supp}\theta_i}G(y-x)e^{-\eta d_x}dx.$$
But, as in (\ref{eq:part2})
$$\vert\int_{x\in\hbox{Supp}\theta_i}G(y-x)e^{-\eta d_x}dx\vert\leq Ce^{-\eta d_y}.$$
In conclusion, (\ref{eq:Greenbis})-(\ref{eq:part7}) give for all $\eta\in]0,1[$ and for a constant $C$ independent of $\beta$, of $R_0$ and of $\ep$
\begin{equation}\label{eq:h1partielle}
|\xi(y) e^{\eta d_y}|\leq C(\| h e^{\eta d_x}\|_{\infty}+{1\over\beta}(\|\nabla\xi e^{\eta d_x}\|_{\infty}+\|\xi e^{\eta d_x}\|_{\infty})+Ce^{-\eta_1 R_0}\|\xi e^{\eta d_x}\|_{\infty}\quad\hbox{for all $y\in \Omega_i$}
\end{equation}
such that $d_y\geq R_0$.\\
Now, for all $y$, we have by (\ref{eq:xiinfini})
$$\vert\xi(y)\vert\leq C\| h\|_{\infty}\leq C\| he^{\eta d_y}\|_{\infty}.$$
Thus
$$\vert\xi(y)e^{d_y}\vert\leq C\| he^{\eta d_y}\|_{\infty}e^{\eta R_0},\quad\hbox{for all $y\in \Omega_i$, such that $d_y\leq R_0$}.$$
 Consequently, we have 
\begin{equation}\label{eq:h1partiellebis} 
|\xi(y) e^{\eta d_y}|\leq C(e^{\eta R_0}\| h e^{\eta d_x}\|_{\infty}+{1\over\beta}(\|\nabla\xi e^{\eta d_x}\|_{\infty}+\|\xi e^{\eta d_x}\|_{\infty})+e^{-\eta_1 R_0}\|\xi e^{\eta d_x}\|_{\infty})
\quad\hbox{for all $y\in \Omega_i$}.
\end{equation}
This estimate being valid for all $i$, we deduce that it is true for every $y\in]-\frac{\pi}{\ep},\frac{\pi}{\ep}[\times\R^{N-1}$.\\

Since $\nabla \xi$ appears in (\ref{eq:h1partiellebis}), we need to prove that 
\begin{equation}\label{eq:h1suite}
|\nabla\xi (y)e^{\eta d_y}|\leq C(e^{\eta R_0}\| h e^{\eta d_x}\|_{\infty}+{1\over\beta}(\|\nabla\xi e^{\eta d_x}\|_{\infty}+\|\xi e^{\eta d_x}\|_{\infty})+e^{-\eta_1 R_0}\|\xi e^{\eta d_x}\|_{\infty})\quad\hbox{}
\end{equation}
for all $y\in ]-\frac{\pi}{\ep},\frac{\pi}{\ep}[\times\R^{N-1}$.\\
Let $y\in\Omega_i$ be given such that $d_y\geq R_0$. Let us denote $u=\xi e^{\eta \vert x-(\frac{a_{\ep}^i}{\ep},0)\vert}$. When $d_x=\vert x-(\frac{a_{\ep}^i}{\ep},0)\vert$, u satisfies the equation
$$-\Delta u+u(1-\eta^2-p\overline u_{\ep}^{p-1})+2\eta\nabla \vert x-(\frac{a_{\ep}^i}{\ep},0)\vert.\nabla u= e^{\eta d_x}h.
$$

 We use Theorem 9.13 of \cite{GT}, with $\Omega'=\{x,\vert x-y\vert\leq {1\over2}\}$ and $\Omega=\{x,\vert x-y\vert< {1}\}$. Without loss of generality, we suppose that $x-(\frac{a_{\ep}^i}{\ep},0)\neq0$ for all $x\in\Omega$. We get
\begin{equation}\label{eq:nabla}
\|\nabla u\|_{L^{\infty}(\Omega')}\leq C(\| u\|_{L^{\infty}(\Omega)}+\|e^{\eta \vert x-(\frac{a_{\ep}^i}{\ep},0)\vert}h\|_{L^{\infty}(\Omega)}).
\end{equation}
We may suppose $y_1>\frac{a_{\ep}^i}{\ep}$. Let $x$ be such that $\vert x-y\vert\leq1$. Then either $d_x=\vert x_1-\frac{a_{\ep}^{i}}{\ep}\vert$ or $d_x=\vert x_1-\frac{a_{\ep}^{i+1}}{\ep}\vert$.
In the latter case, we have
$$\vert\vert x_1-\frac{a_{\ep}^{i+1}}{\ep}\vert-\vert x_1-\frac{a_{\ep}^{i}}{\ep}\vert\vert\leq 1.$$
Consequently, there exists $C$ independent of $\ep$ such that
$$\hbox{$\forall x\in\Omega$}\quad 
\vert d_x-\vert x-(\frac{a_{\ep}^i}{\ep},0)\vert\vert\leq C 
.$$
Thus, (\ref{eq:nabla}) gives
\begin{equation}\label{eq:h1fin}
\vert\nabla(\xi e^{\eta d_y})(y)\vert\leq C(\|\xi e^{\eta d_x}\|_{L^{\infty}(\Omega)}+\| h e^{\eta d_x}\|_{L^{\infty}(\Omega)}).
\end{equation}
Then (\ref{eq:h1partiellebis}) and (\ref{eq:h1fin}) lead to (\ref{eq:h1suite}). 
Now (\ref{eq:h1suite}) and (\ref{eq:h1partielle})  (which is valid for all $y$)  give
\begin{equation}\label{eq:hj2}
(\vert\xi(y)\vert+\vert\nabla\xi(y)\vert )e^{\eta d_y}\leq C(e^{\eta R_0}\|he^{d_x}\|_{\infty}+ {1\over\beta}(\|\xi e^{d_x}\|_{\infty}+\|\nabla\xi e^{d_x}\|_{\infty})+e^{-\eta_1 R_0}\|\xi e^{\eta d_x}\|_{\infty}\quad\forall y.
\end{equation}

Now we choose $\beta$ and $R_0$ large enough to get (\ref{eq:hj1}). We have proved the proposition.\\

\begin{proposition}\label{barrierefonctionpropre}
Let $\varphi$ be an eigenfunction of $\LL$, associated with an eigenvalue $\lambda$ which tends to $\overline\lambda\neq1$. Let us suppose that $\|\varphi\|_{L^{\infty}(]-\frac{\pi}{\ep},\frac{\pi}{\ep}[\times\R^{N-1})}=1$.
Then for all $\eta\in]0,1[$ there exists $C>0$, independent of $\ep$, such that for $\ep$ small enough (depending on $\overline\lambda$)
\begin{equation}\label{eq:fp}
\vert\varphi(x)\vert+\vert\nabla\varphi(x)\vert\leq {C}e^{-\eta\hbox{dist}(x,\cup_{i=0}^{k+1}\{(\frac{a_{\ep}^i}{\ep},0)\})},
\end{equation}
where we use the notation : $a_{\ep}^0=a_{\ep}^k-2\pi$ and $a_{\ep}^{k+1}=a_{\ep}^1+2\pi$.\\
Moreover
\begin{equation}\label{eq:normes}
C_1\leq\|\varphi\|_{H^1}\leq C_2
\end{equation}
where $C_1$ and $C_2$ are some positive real numbers independent of $\ep$.\\
Let $\xi$ be defined in (\ref{eq:theo}). Then
\begin{equation}\label{eq:xip}
\vert\xi(x)\vert+\vert\nabla\xi(x)\vert\leq C(e^{-\underline\sigma}+\sum_{j=1}^k\vert\lambda_j\vert )e^{-\eta\hbox{dist}(x,\cup_{i=0}^{k+1}\{(\frac{a_{\ep}^i}{\ep},0)\})}.
\end{equation}
\end{proposition}

{\bf Proof}  To prove (\ref{eq:fp}), we follow the proof of (\ref{eq:hj1}) in Proposition \ref{appendix1}, in which we let $h=0$ and we replace $\overline u_{\ep}^{p-1}$ by ${\overline u_{\ep}^{p-1}/ 1-\lambda}$.\\
We find, for $y\in \Omega_i$ such that $d_y\geq R_0$
$$\vert\varphi(y)e^{\eta d_y}\vert+\vert\nabla\varphi(y)e^{\eta d_y}\vert\leq{C\over\beta}(\|\varphi e^{d_x}\|_{\infty}+\|\nabla\varphi e^{d_x}\|_{\infty})+Ce^{-\eta_1 R_0}\|\varphi e^{\eta d_x}\|_{\infty}
$$
while for $y\in \Omega_i$ such that $d_y\leq R_0$
$$\vert\varphi(y)e^{\eta d_y}\vert+\vert\nabla\varphi(y)e^{\eta d_y}\vert\leq Ce^{\eta R_0}$$
where the constant $C$ is independent of $\beta$ and of $R_0$ but are dependent on $\lambda$.
We choose $R_0$ and $\beta$ large enough to obtain (\ref{eq:fp}).\\

Let us prove (\ref{eq:normes}).\\
We have 
 $$(1-\lambda)\|\varphi\|^2_{H^1}=p\int_{\frac{S^1}{\ep}\times\R^{N-1}}\overline u_{\ep}^{p-1}\varphi^2dx.$$
 In view of (\ref{eq:fp}) and of (\ref{eq:norme}), we may use the Lebesgue Theorem to obtain that
 $$\|\varphi\|_{H^1}\not\rightarrow0,\quad i=0,...,k.$$
Now we have
$$\LL\xi=\LL\frac{\partial v_i}{ \partial x_1}-\sum_{j=1}^kc_j\lambda_j(-\Delta\varphi_j+\varphi_j)=\sum_{l\in\Z}p(\overline u_{\ep}^{p-1}-U^{p-1}_{i,l})\frac{\partial U_{i,l}}{\partial x_1}+p\sum_{j=1}^kc_j{\lambda_j\over 1-\lambda_j}\overline u_{\ep}^{p-1}\varphi_j.
$$
If we write $\LL\xi=h$, then, for all $\eta\in]0,1[$ $$\vert h(x)\vert \leq C_{\eta} (e^{-\underline\sigma}+\sum_{j=1}^k\vert\lambda_j\vert )e^{-\eta d_x}.$$
We use Proposition \ref{appendix1} to obtain (\ref{eq:xip}).\\
We have proved the proposition.\\

\begin{proposition}\label{convolution1} Let $C$ and $A$ be positive real numbers. 
Let $f$ and $g$ be functions which verify the following property, for $\vert x\vert>A$
$$\vert f(x)\vert\leq Ce^{-\vert x\vert}\quad, \quad\vert g(x)\vert\leq Ce^{-\vert x\vert}\vert x\vert^{1-N\over2}.$$
Let $a>b>0$. Let $y_0$ be such that $\vert y_0\vert\rightarrow+\infty$ and $\alpha={\vert y_0\vert\over2}$. Then
\begin{equation}
\vert\int_{\alpha^{1\over2}<\vert x\vert<\alpha}f^a(x)g^b(x-(y_0,0))dx\vert=o(\vert y_0\vert^{b\frac{1-N}{2}}e^{-b\vert y_0\vert}).
\end{equation}
\end{proposition}

{\bf Proof}  We easily see that if $\vert x\vert<\alpha$,
$$\vert x-(y_0 ,0)\vert=\vert y_0\vert\vert{x\over \vert y_0\vert}-(1,0)\vert\geq \frac{1}{2}\vert y_0\vert.$$
Thus, assuming $y_0>0$
$$\vert\int_{\alpha^{1\over2}<\vert x\vert<\alpha}f^a(x)g^b(x-(y_0,0))dx\vert\leq C \vert y_0\vert^{b\frac{1-N}{2}}
\int_{\alpha^{1\over2}<\vert x\vert<\alpha}e^{-a\vert x\vert-b\vert x-(y_0,0)\vert}dx$$
$$\leq C\vert y_0\vert^{b\frac{1-N}{2}}\int_{\alpha^{1\over2}<r<\alpha}r^{N-1}\int_{S^{N-1}} e^{-ar}e^{-b\sqrt{(rz_1-y_0)^2+r^2\sum_{i=2}^N z_i^2}}dr d\mu(z)$$
$$\leq C\vert y_0\vert^{b\frac{1-N}{2}}\int_{\alpha^{1\over2}<r<\alpha}r^{N-1}\int_{S^{N-1}} e^{-ar}e^{-b\sqrt{r^2-2rz_1y_0+y_0^2}}dr d\mu(z)$$
$$\leq C\vert y_0\vert^{b\frac{1-N}{2}}\int_{\alpha^{1\over2}<r<\alpha}r^{N-1}\int_{S^{N-1}} e^{-ar}e^{-b\vert y_0-r\vert}dr d\mu(z)$$
$$\leq C\vert y_0\vert^{b\frac{1-N}{2}}e^{-b\vert y_0\vert}e^{(b-a)\alpha^{1\over2}}\alpha^N$$
and $e^{(b-a)\alpha^{1\over2}}\alpha^N\rightarrow0$, since $a>b$.\\

\begin{proposition}\label{convolution2}
Let $f$ and $g$ be smooth functions and $C$, $C_1$, $C_2$ and $A$ be positive real numbers which verify, for $\vert x\vert>A$
$$0\leq f(x)\leq Ce^{-\vert x\vert}\quad;\quad 
C_1e^{-\vert x\vert}\vert x\vert^{1-N\over2}\leq  g(x)\leq C_2e^{-\vert x\vert}\vert x\vert^{1-N\over2}.$$
Let us define
$$\Omega_{i}=\{x\in[-{\pi\over\ep},{\pi\over\ep}]\times\R^{N-1};\hbox{dist}( x,\cup_{l=0}^{k+1}\{(\frac{a_{\ep}^l}{\ep},0)\}=\vert x-(\frac{a_{\ep}^i}{\ep},0)\vert\}.$$
Let $a>b>0$. Let $i\neq j$ and let $y_0={a_{\ep}^j-a_{\ep}^i\over\ep}$. \\
If $f(x)\geq0$ and $f\neq0$, then there exist positive real numbers $C_1$ and $C_2$ such that, for $i=0,...,k+1$
\begin{equation}\label{eq:convol1}
(1+o(1))C_1^bC_0 e^{-b\vert y_0\vert}\vert y_0\vert^{b{1-N\over2}}\leq\int_{\Omega_{i}}f^a(x-(\frac{a_{\ep}^i}{\ep},0))g^b(x-(\frac{a_{\ep}^j}{\ep},0))dx\leq C_2^b C_0e^{-b\vert y_0\vert}\vert y_0\vert^{b{1-N\over2}}(1+o(1)),
\end{equation}
where
$$C_0=\int_{\R^N}f^a(x)dx.$$
The estimate (\ref{eq:convol1}) holds true if we replace $\Omega_i$ by the set
$$\Omega_i^+=\{x\in\Omega_i;\quad x_1>\frac{a_{\ep}^i}{\ep}\}$$
while $C_0$ is replaced by
$$\int_{\R^N, x_1>0}f^a(x)dx.$$
\end{proposition}

{\bf Proof.} Let $$\alpha= \frac{\vert y_0\vert}{2}.$$

For $x$ such that $\vert x-(\frac{a_{\ep}^i}{\ep},0)\vert<\alpha^{1\over2}$, we have
$$\vert y_0\vert -\vert x-(\frac{a_{\ep}^i}{\ep},0)\vert\leq \vert x-(\frac{a_{\ep}^j}{\ep},0)\vert\leq\vert y_0\vert +\vert x-(\frac{a_{\ep}^i}{\ep},0)\vert,$$
thus
$$\vert x-(\frac{a_{\ep}^j}{\ep},0)\vert= (1+o(1))\vert y_0\vert.$$ 
We write
$$\int_{\vert x-(\frac{a_{\ep}^i}{\ep},0)\vert<\alpha^{1\over2}}f^a(x-(\frac{a_{\ep}^i}{\ep},0))g^b(x-(\frac{a_{\ep}^j}{\ep},0))dx$$
$$\geq C_1^b(1+o(1))\vert y_0\vert^{b{(1-N)\over2}}e^{-b\vert y_0\vert}\int_{\vert x-(\frac{a_{\ep}^i}{\ep},0)\vert<\alpha^{1\over2}}f^a(x-(\frac{a_{\ep}^i}{\ep},0))dx
$$
$$\geq C_1^bC_0(1+o(1))\vert y_0\vert^{b{(1-N)\over2}}e^{-b\vert y_0\vert}$$
and
$$\int_{\vert x-(\frac{a_{\ep}^i}{\ep},0)\vert<\alpha^{1\over2}}f^a(x-(\frac{a_{\ep}^i}{\ep},0))g^b(x-(\frac{a_{\ep}^j}{\ep},0))dx$$
$$\leq C_2^b(1+o(1))\vert y_0\vert^{b{(1-N)\over2}}e^{-b\vert y_0\vert}\int_{\vert x-(\frac{a_{\ep}^i}{\ep},0)\vert<\alpha^{1\over2}}f^a(x-(\frac{a_{\ep}^i}{\ep},0))dx
$$
$$\leq C_2^bC_0(1+o(1))\vert y_0\vert^{b{(1-N)\over2}}e^{-b\vert y_0\vert}.$$
Moreover, in view of Proposition \ref{convolution1}, we have
$$\int_{x\in\Omega_i;\alpha^{1\over2}\leq\vert x-(\frac{a_{\ep}^i}{\ep},0)\vert<\alpha}f^a(x-(\frac{a_{\ep}^i}{\ep},0))g^b(x-(\frac{a_{\ep}^j}{\ep},0))dx=o(1)\vert y_0\vert^{b{(1-N)\over2}}e^{-b\vert y_0\vert}.$$
Last, if $x\in \Omega_i$, we have $\vert x-(\frac{a_{\ep}^j}{\ep},0)\vert\geq \alpha$. We deduce that
$$\int_{x\in\Omega_i; \vert x-(\frac{a_{\ep}^i}{\ep},0)\vert>\alpha}f^a(x-(\frac{a_{\ep}^i}{\ep},0))g^b(x-(\frac{a_{\ep}^j}{\ep},0))dx$$
$$\leq C_2^bC^ae^{-b\vert y_0\vert}\vert y_0\vert^{b{1-N\over2}}\int_{x\in\Omega_i; \vert x-(\frac{a_{\ep}^i}{\ep},0)\vert>\alpha}e^{b\vert x-(\frac{a_{\ep}^i}{\ep},0)\vert}f^a(x-(\frac{a_{\ep}^i}{\ep},0))dx$$
$$=o(e^{-b\vert y_0\vert}\vert y_0\vert^{b{1-N\over2}}),$$
since $a>b$.\\
We have proved the proposition.\\

\begin{corollary}\label{equivalence}
If $$\lim_{\vert x\vert\rightarrow+\infty} g(x)e^{\vert x\vert}\vert x\vert^{N-1\over2}= L$$
then $$\lim_{\ep\rightarrow0} \int_{\Omega_i}f^a(x-(\frac{a_{\ep}^i}{\ep},0))g^b(x-(\frac{a_{\ep}^j}{\ep},0))dx e^{b\vert y_0\vert}\vert y_0\vert^{b{N-1\over2}}=LC_0$$
where $C_0$ is defined in Proposition \ref{convolution2}. Moreover, we can replace $\Omega_i$ by $\Omega_i^+$.
\end{corollary}

 \begin{lemma}\label{taylor} Let $k$ be a positive integer. Let $a$ and $b$ be real numbers, with $a>0$.
If $p>1$ is given and $k=[p]$. There exists $C$ independent of $a$ and $
b$ such that
$$\vert -(a+b)_+^p+a^p+...+{p....(p-k+1)\over k!}a^{p-k}b^{k}\vert\leq C\vert b\vert^p.$$
\end{lemma}
{\bf Proof}  First, let us suppose that $b<0$. If $a+2b<0$, then $a<2\vert b\vert$, and the claim is true.\\
If $a+2b>0$, then $a+b>0$ and we write
\begin{equation}\label{eq:taylor2}
\vert -(a+b)^p+a^p+...+{p....(p-k+1)\over k!}a^{p-k}b^{k}\vert={p....(p-k)\over (k+1)!}(a+\alpha b)^{p-k-1}b^{k+1}
\end{equation}
where $\alpha\in]0,1[$. But $a+\alpha b>-2 b+\alpha b=(2-\alpha)\vert b\vert$, thus 
$$\vert -(a+b)^p+a^p+...+{p....(p-k+1)\over k!}a^{p-k}b^{k}\vert\leq{p....(p-k)\over (k+1)!}(2-\alpha )^{p-k-1}\vert b\vert^{p}$$
that gives the claim.\\
Secondly, let us suppose that $b>0$. If $a<b$, the claim is true. If $a>b$, we use (\ref{eq:taylor2}) again and we obtain the claim.\\


\end{document}